\documentclass[12pt]{article}
\usepackage{e-jc}

\newcommand{\Z}{{\bf Z}}
\newcommand{\R}{{\bf R}}

\newcommand{\cH}{{\cal H}}
\newcommand{\cM}{{\cal M}}
\newcommand{\cR}{{\cal R}}

\newcommand{\cT}{{\cal T}}

\newcommand{\meet}{{\wedge}}
\newcommand{\join}{{\vee}}

\newcommand{\symmdiff}{{\bigtriangleup}}

\usepackage{amsmath,graphicx}

\dateline{Jan 1, 2020}{Jan 2, 2020}{TBD}


\Copyright{The author. Released under the CC BY-ND license (International 4.0).}

\title{Lattice structure for orientations of graphs}

\author{James Propp\authornote{1} }
\authortext{1}{Department of Mathematics and Statistics, University of Massachusetts Lowell, Lowell, Massachusetts, USA (\email{jamespropp@gmail.com}).}

\begin{document}

\maketitle


\begin{abstract}
In 1986, Oliver Pretzel studied
the set of orientations of a connected finite graph $G$
and showed that any two such orientations having the same
flow-difference around all closed loops can be obtained
from one another by a succession of local moves of a simple type.
Here I show that the set of orientations of $G$
having the same flow-differences around all closed loops
can be given the structure of a distributive lattice.
When the graph is drawn on the plane, a dual version of the construction
puts a distributive lattice structure on the set of orientations 
of $G$ having the same indegrees at all vertices.
In both settings, adjacent lattice-elements are related
by simple local moves.
This construction unifies earlier, similar constructions
in combinatorics and statistical mechanics.
It also gives rise to an interesting lattice structure
on spanning trees.
This article is an updated version of a preprint
originally distributed in 1993.

\end{abstract}

\section{Statement of results.}

This paper extends Thurston's definition of height functions,
originally limited to domino and lozenge tilings~\cite{Thur}, 
in a number of different directions.
These extensions can be treated in a unified fashion by using a framework 
introduced by Oliver Pretzel~\cite{Pr1} a decade earlier. 

The main result (Theorem 1)
asserts that the set of orientations of a graph
having the same flow-differences around all closed loops
can be given the structure of a distributive lattice.
A full statement is given later in this section,
and a proof is given in Section 2.

In the case where the graph is drawn on the plane,
we can take the dual of the construction, 
obtaining a result in which the circulation around a loop 
is replaced by the net outdegree at a vertex
(the discrete analogue of the duality between curl and divergence).
A full statement is given later in this section,
and a proof is given in Section 3.

A still more specialized result gives a lattice structure
on the set of spanning trees of a finite connected graph 
drawn on the sphere.
A full statement is given later in this section,
and a proof is given in Section 4.

Finally, Section 5 provides some historical context.
In the decades in between the original 
issuance of this article in preprint form and its publication,
many developments arose building on these ideas,
some of which were rediscovered by other researchers
who were unaware of my preprint.
I, in turn, had been unaware of related work that preceded my own.

I now proceed to introduce terminology and conventions
that will enable me to state the three main theorems.

All graphs are assumed to be finite and connected.
Multiple edges are allowed, but self-loops are not permitted.
Let $G$ be a graph with vertex-set $V$ and edge-set $E$.
An {\bf orientation} of the edge $e$ is a triple $(e,v,w)$
where $v$ and $w$ are the endpoints of $e$;
if $\vec{e} = (e,v,w)$, then we say $\vec{e}$ points from $v$ to $w$,
and we define the reversal of $\vec{e}$ as $-\vec{e} = (e,w,v)$.
Let $\vec{E}$ be the set of directed edges of $G$,
i.e., the set containing all $2|E|$ orientations of edges of $G$.
Writing $E = \{e_1,e_2,...\}$,
we define an {\bf orientation} $R$ of $G$ as a set of directed edges
$\{\vec{e}_1,\vec{e}_2,...\} \subset \vec{E}$
where each $\vec{e}_i$ is an orientation of the edge $e_i$.

When the graph $G$ is planar,
we can either fix an embedding on the sphere with a marked face $f^*$
or fix an embedding in the plane 
(where the unbounded face corresponds to $f^*$).
Some arguments are easier to see on the sphere
(some do not even require attending to $f^*$)
and others are easier to see in the plane,
so I will freely (but not, I hope, capriciously)
pass between the two sorts of pictures.
I will say that $G$ is ``drawn on the plane''
(equivalently, ``is a plane graph'')
or ``drawn on the sphere'' to indicate which viewpoint
I think is helpful, but in both cases it should be understood
that there are no crossing edges (except at the end of Section 4).

A {\bf directed path} in $G$
with initial vertex $v_0$ and terminal vertex $v_n$
is a sequence of directed edges
$(e_1,v_0,v_1),(e_2,v_1,v_2),...,(e_n,v_{n-1},v_n)$.
The sequence of edges
$e_1,e_2,...,e_n$ that is associated with a directed path
is called simply a {\bf path}.
A path with no repeated vertices is {\bf simple}.
A directed path whose initial vertex and terminal vertex coincide
is a {\bf directed cycle}; the edges associated with it form a {\bf cycle}.
A cycle is {\bf simple}
if its only repeated vertex is the initial/terminal vertex.
We say a directed edge is a {\bf forward edge}
(relative to the orientation $R$) if it belongs to $R$,
and a {\bf backward edge} if it does not;
a path composed entirely of forward edges is a {\bf forward path}.
Given vertices $v,w$ of $G$, we say that $w$ is {\bf accessible} from $v$
(relative to $R$) if $R$ contains a forward path from $v$ to $w$.
Mutual accessibility is an equivalence relation;
we call its equivalence classes {\bf accessibility classes}.
If all equivalence classes are of size 1,
we say that $R$ is {\bf acyclic}.
In most cases of interest, $R$ is acyclic,
but we will treat the general case.

Let $R$ be an orientation of $G$.
An accessibility class $A$ will be called {\bf maximal}
if all directed edges in $R$ between $A$ and its complement $A^c$
point towards $A$, and {\bf minimal} if they all point toward $A^c$.
If $A$ is maximal,
the process of reversing the orientations of all the directed edges
between $A$ and $A^c$ (thereby making $A$ minimal rather than maximal)
is called {\bf pushing down} $A$ (see~\cite{Pr1,Pr2}).
Figure 1 gives two examples of pushing down
that are inverse to one another.
In the first example, we push down a single vertex $v_1$;
in the second, we push down a non-trivial accessibility class
that is the complement of $\{v_1\}$, namely $\{v_2,v_3,v_4\}$.
Pretzel calls the former sort of move ``pushing down''
and the latter sort of move ``PDCE'',
but I prefer to use the suggestive term ``pushing down'' for both.
Vertices in the accessibility class 
that is about to be pushed down are circled.
The idea of pushing down a maximal vertex appears to be due to Mosesian 
(see the three articles cited in 
the bibliographies of~\cite{Pr1} and~\cite{Pr3}).

\begin{figure}
\begin{center}
\includegraphics[width=4in]{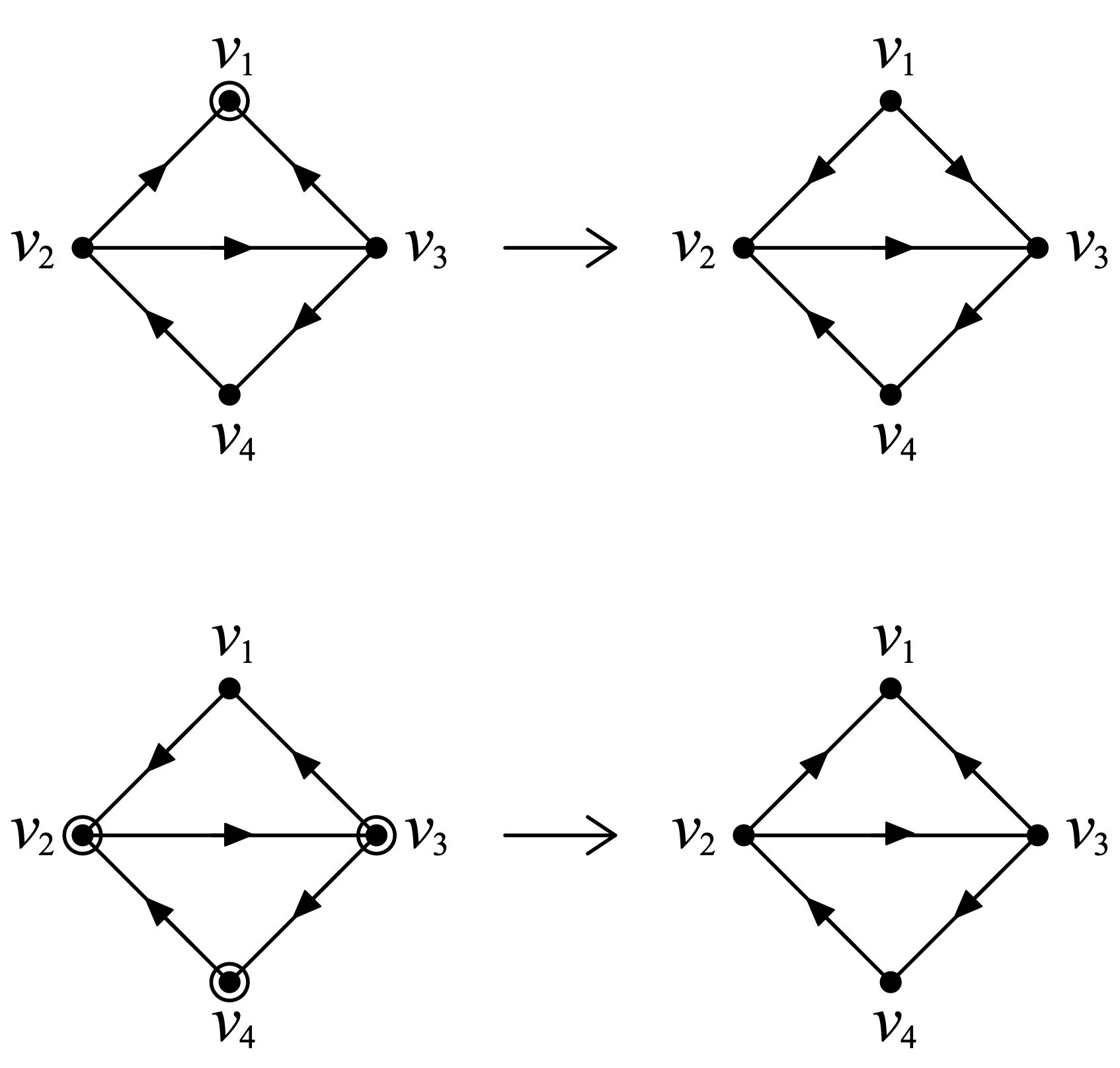}
\end{center}
\caption{Two examples of pushing down.} 
\end{figure}

If $C$ is a directed cycle in $G$,
define the {\bf circulation} of $R$ around $C$
(Pretzel calls this the {\bf flow-difference} around the cycle)
as the integer $|C_R^+|-|C_R^-|$,
where $C_R^+$ is the set of forward edges of $C$
and $C_R^-$ is the set of backward edges of $C$.
Define the circulation of $R$ as the function $c=c_R$
that associates to each directed cycle $C$
the circulation $c(C)$ of $R$ around $C$;
say also that $R$ is a
\boldmath $c${\bf -orientation}.
\unboldmath
(One can see circulation as a discrete analogue of
the scalar curl of a two-dimensional vector field,
much as the number of outward-pointing edges at a vertex
minus the number of inward-pointing edges at the vertex
can be seen as a discrete analogue of divergence.)
Note that the problem of determining
whether some orientation $R$ is a $c$-orientation
reduces to the problem of evaluating the circulation of $R$
around all cycles $C$ belonging to a cycle-basis.
If $R$ is a $c$-orientation,
then $|C_R^+| + |C_R^-| = |C|$ and $|C_R^+| - |C_R^-| = c(C)$,
so $|C_R^+|=\frac12(|C|+c(C))$ and $|C_R^-|=\frac12(|C|-c(C))$.
We say that a circulation $c$ is {\bf feasible}
if there exists at least one $c$-orientation of $G$.

It is easy to check that pushing down does not affect
the circulation of an orientation of $G$.
Indeed, let $R$ be an orientation with
a maximal accessibility class $A$,
and let $C$ be any directed cycle in $G$.
Since $C$ must contain equal numbers of directed edges
going from $A$ to $A^c$ and going from $A^c$ to $A$,
reversing the orientations of all the directed edges
joining $A$ and $A^c$ has no effect
on the circulation of $R$ around $C$.

I now state the main result of this article (proved in Section 2).

\begin{theorem} \label{thm1}
Let $\cR$ be the (non-empty) set of orientations
of a finite connected graph $G$
that have a fixed circulation $c$,
and let $A^*$ be an accessibility class of $G$.
If we say that one $c$-orientation $R$ covers another $c$-orientation $S$
exactly when $S$ is obtained from $R$ by pushing down
at a maximal accessibility class other than $A^*$,
then the covering relation makes $\cR$ into a distributive lattice.
\end{theorem}

\begin{figure}
\begin{center}
\includegraphics[width=3in]{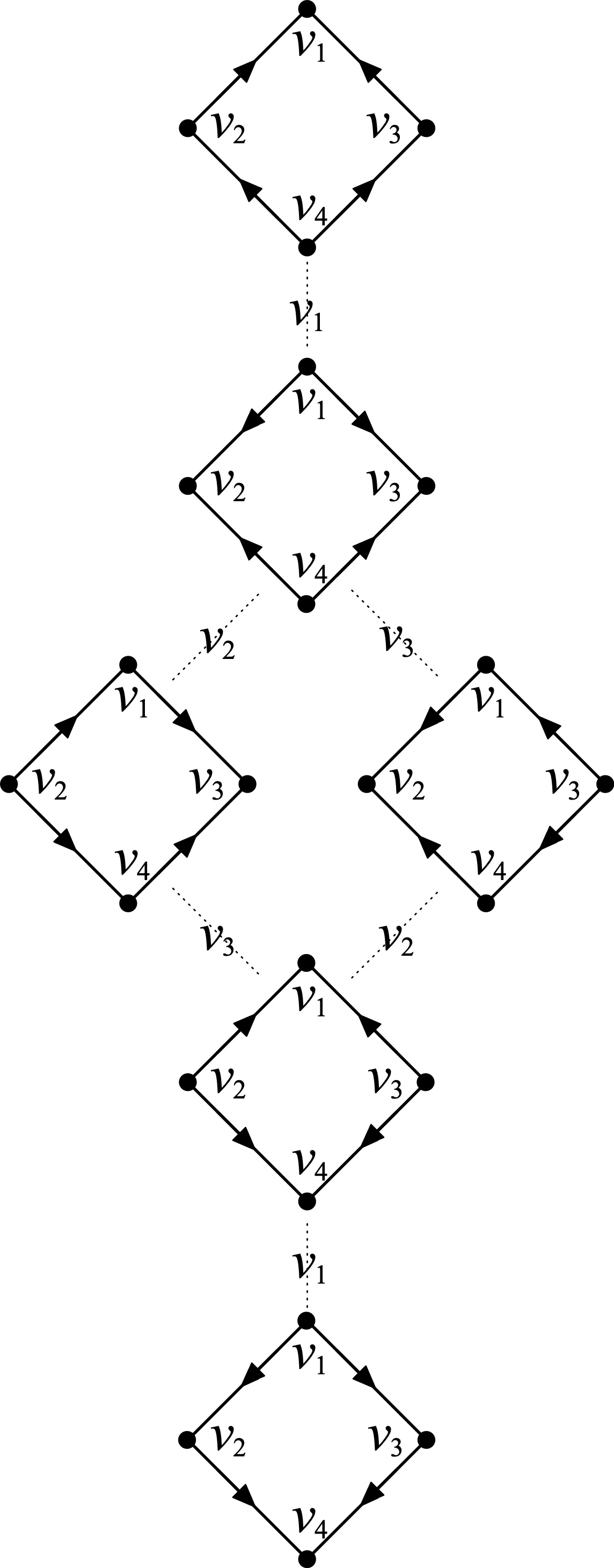}
\end{center}
\caption{Orientations of a rooted 4-cycle.}
\end{figure}

Figure 2 gives an example of such a distributive lattice $\cR$
in the case where $G$ is a cycle of length 4,
$\cR$ is the set of orientations of $G$
in which the circulation around the cycle is zero,
and $A^* = \{v_4\}$ (the ``root'').
The figure shows the Hasse diagram of the lattice.
Each edge (shown as a dotted line) corresponds to a pushing down move;
for instance, the orientation shown at the top of the figure
turns into the orientation shown just below it
if the sink at $v_1$ is pushed down. 
This example is discussed in greater detail later (Example \ref{E11}).

The proof of Theorem~\ref{thm1} will introduce height-functions
(inspired by and generalizing the height-function
introduced by Thurston~\cite{Thur}) to aid the analysis.
Using height-functions, I will prove that
if $R$ and $S$ are $c$-orientations of $G$,
then $S$ can be obtained from $R$ by
performing at most $N(N-1)/2$ pushing-down operations,
where $N$ is the number of vertices of $G$.
This result strengthens the theorem of Pretzel~\cite{Pr1},
who showed that a finite number of pushing-down operations suffices.
The bound $N(N-1)/2$ is best possible; 
see Propositions \ref{PK} and \ref{PL}.

Now suppose $G$ is a connected bipartite graph drawn on the sphere
and $d$ is a function from the vertex-set of $G$
to the non-negative integers.
A \boldmath $d${\bf -factor} \unboldmath
of $G$ is a subgraph of $G$
in which each vertex $v$ has degree $d(v)$.
For example, if $d(v)=1$ for all vertices $v$,
then a $d$-factor of $G$
is just a $1$-factor, or (perfect) matching, of $G$.
This being the most important case,
I will use the letter $M$ (for matching)
to denote a $d$-factor.
We formally regard $M$ as a set of edges;
thus $M^c$ denotes the set of edges of $G$ not in $M$.

Note that the edges of $G$
divide the sphere into simply-connected regions, or {\bf faces},
with an even number of sides,
if we agree to double-count edges 
that are internal to a face, as in Figure 3.
We will pick some particular face $f^*$
that will play a special role in what follows.
In depicting $G$ in the plane (on the page),
we will have the unbounded, external face correspond to $f^*$.
We will give each face a preferred orientation,
which is counterclockwise (on the sphere);
note that when one draws the graph on the plane,
the preferred orientation of the outer face $f^*$
appears to be clockwise.

\begin{figure}
\begin{center}
\includegraphics[width=2in]{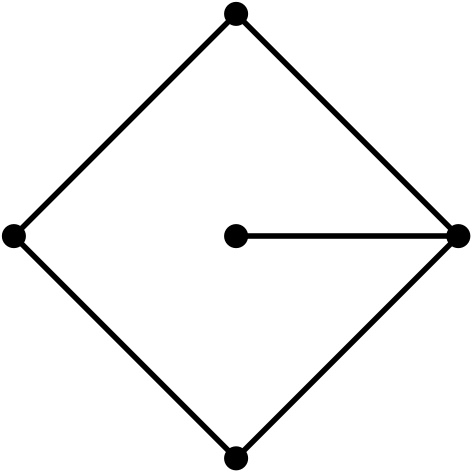}
\end{center}
\caption{A six-sided face.}
\end{figure}

Assume that every edge of $G$
belongs to some $d$-factors but not to others.
Note that if this is not true,
the $d$-factors of $G$
are in 1-to-1 correspondence with
the $d'$-factors of some subgraph $G'$,
where every edge of $G'$ belongs to some $d'$-factors but not to others.
Thus, we are not losing any generality
in making the assumption.

Color the vertices of $G$ alternately black and white.
An {\bf elementary cycle} in $G$
is a simple cycle that encircles a single face.
An {\bf alternating cycle} in $G$ relative to a $d$-factor $M$
is an elementary cycle in $G$
in which the edges alternately belong to $M$ and $M^c$.
Call the cycle {\bf positive}
if the edges in the cycle,
when directed from black vertices to white vertices,
circle the face in a counterclockwise direction,
and {\bf negative}
if they circle in a clockwise direction.
(A cycle need not be positive or negative.)
A {\bf face twist}
is the operation of removing from a $d$-factor
some edges that form an alternating cycle around a face
and inserting the complementary edges.
More specifically, {\bf twisting down}
is the operation that converts a positive alternating cycle
into a negative alternating cycle around the same face,
and {\bf twisting up} is its inverse.
Note that twisting converts a $d$-factor into another $d$-factor.
Figure 4 shows the three 1-factors of a particular graph,
in which the faces having positive or negative alternating cycles
are marked $+$ or $-$, respectively;
an arrow from a 1-factor $M$ to a 1-factor $N$
expresses the fact that $N$ can be obtained from $M$ by twisting down
a positive alternating cycle.

\begin{figure}
\begin{center}
\includegraphics[width=0.6\textwidth]{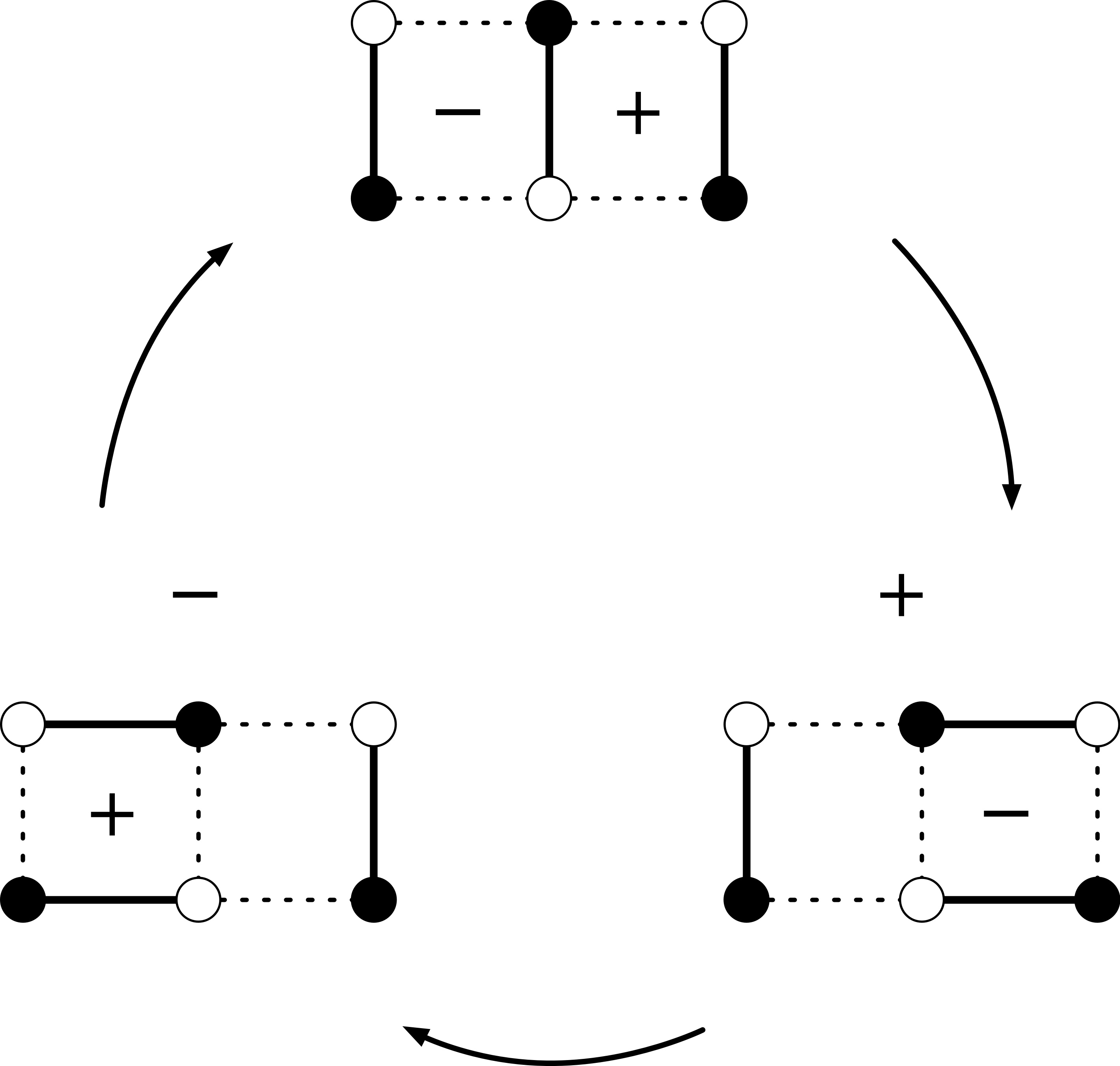}
\end{center}
\caption{Twisting down a 1-factor.}
\end{figure}

The following consequence of Theorem \ref{thm1},
proved in Section 3,
gives us a lattice structure for the set of $d$-factors of $G$.

\begin{theorem} \label{thm2}
Let $\cM$ be the (non-empty) set of $d$-factors
of a finite graph $G$ drawn on the sphere.
If we say that one $d$-factor $M$ covers another $d$-factor $N$
exactly when $N$ is obtained from $M$ by twisting down
at a face other than $f^*$,
then the covering relation makes $\cM$ into a distributive lattice.
\end{theorem} 

Lastly, let $G$ be a connected graph 
drawn on the sphere with a special face $f^*$.
Let $v^*$ be a vertex of $G$ on the boundary of the face $f^*$.
A {\bf spanning tree} of $G$ is a collection $T$ of edges of $G$
such that, for any two vertices $v,w$ of $G$,
there is a unique simple path between $v$ and $w$ using edges of $T$.

Let $v$ be a vertex of $G$ and $e$ an edge containing $v$.
If we imagine an ant that travels in a small clockwise circle centered at $v$,
we see that after it crosses $e$, it remains in some face $f$ of $G$
before crossing the next edge $e'$ of $G$.
We say that $e'$ is the {\bf clockwise successor} of $e$ at $v$.
Fix $v,e,f,e'$ as described, and suppose $e \neq e'$.
If $T$ is a spanning tree of $G$,
we will say that the angle $eve'$ is {\bf positively pivotal}
if the following conditions are satisfied:
\begin{enumerate}
\item[(1)] $e \in T$ and $e' \not\in T$;
\item[(2)] the symmetric difference $T' = T \symmdiff \{e,e'\}$
is a spanning tree of $G$;
\item[(3)] the simple path from $v$ to $v^*$ in $T$ contains $e$;
\item[(4)] the simple path from $v$ to $v^*$ in $T'$ contains $e'$; and
\item[(5)] the (unique) simple cycle in $T \cup T'$ separates $f$ from $f^*$.
\end{enumerate}
When $eve'$ is positively pivotal,
the operation $T \mapsto T'$ is called {\bf swinging down} at $v$
through the angle $eve'$.
Note that condition (5) necessitates $f \neq f^*$.
(We may define the notions ``negatively pivotal'' and
``swinging up'' in an analogous way.)

\begin{figure}
\begin{center}
\includegraphics[width=4in]{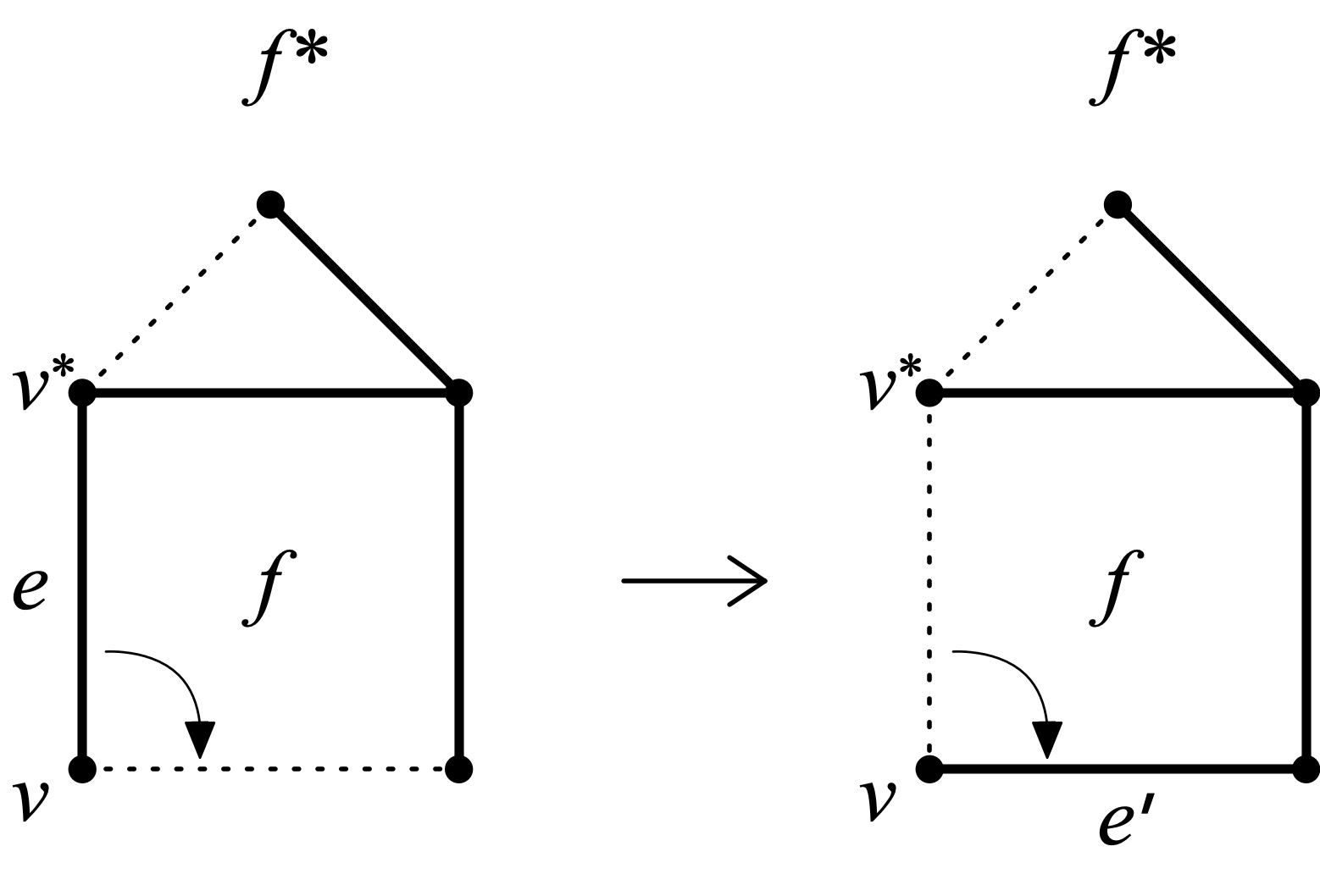}
\end{center}
\caption{A swinging-down move.}
\end{figure}

Figure 5 shows an example of a swinging-down move.
The special vertex $v^*$ is marked with a large black dot;
the swinging-down move takes place at the lower-left angle
of the four-sided face.

The following consequence of Theorem \ref{thm1},
proved in Section 4,
gives us a lattice structure for the set of spanning trees of $G$.

\begin{theorem} \label{thm3}
Let $\cT$ be the set of spanning trees
of a finite connected graph $G$ drawn on the sphere,
with $f^*$ a face of $G$
and $v^*$ a vertex of $G$ incident with $f^*$.
If we say that one tree $T \in \cT$ covers another tree $T' \in \cT$
exactly when $T'$ is obtained from $T$ by swinging down,
then the covering relation makes $\cT$ into
a distributive lattice.
\end{theorem}

\section{Orientations of graphs.}

For the reader's convenience, I restate Theorem \ref{thm1}.

\bigskip

\noindent
{\bf Theorem \ref{thm1}:} {\em Let $\cR$ be the (non-empty) set of orientations
of a finite connected graph $G$
that have a fixed circulation $c$,
and let $A^*$ be an accessibility class of $G$.
If we say that one $c$-orientation $R$ covers another $c$-orientation $S$
exactly when $S$ is obtained from $R$ by pushing down
at a maximal accessibility class other than $A^*$,
then the covering relation makes $\cR$ into a poset
that is a distributive lattice.}

\bigskip

I will prove this theorem in a slightly roundabout way.
Specifically, I will use height-functions
to define a lattice structure on $\cR$,
and then verify that the covering relation for this ordering
is the relation given by pushing down.

Note that two vertices $v,w$ of $G$ are mutually accessible
(relative to an orientation $R$) if and only if
there exists a directed cycle of length $n$ containing $v$ and $w$
around which the circulation is $n$.
If $S$ is another $c$-orientation of the graph,
then the circulation of $S$ around this cycle must also be $n$,
implying that the cycle is a forward cycle
and that $v$ and $w$ are mutually accessible relative to $S$.
It follows that any two $c$-orientations of $G$
induce the same partition of the vertex-set $V$ into accessibility classes.
In particular, we can talk about the partition of $V$ induced by
the circulation $c$, and we can describe $c$ as being acyclic (or not).

Given a feasible circulation $c$ of $G$,
let us describe a directed edge $\vec{e}$ as being {\bf forced}
if it belongs to every $c$-orientation of $G$,
and {\bf forbidden} if it belongs to no $c$-orientation of $G$.
Clearly $\vec{e}$ is forced if and only if $-\vec{e}$ is forbidden.

\begin{proposition} \label{PZ}
Let $c$ be a feasible circulation on $G$.
A directed edge $\vec{e}$ is forced or forbidden
if and only if
its endpoints are in the same accessibility class.
\end{proposition}

\begin{proof}
Suppose $\vec{e}=(e,v,w)$
has endpoints $v,w$ that are in the same accessibility class.
Let $R$ be some $c$-orientation of $G$.
If $R$ contains $\vec{e}$,
then, since $R$ contains a path from $w$ to $v$,
$R$ contains a forward cycle that includes the directed edge $\vec{e}$;
all the edges in this cycle, including $\vec{e}$, must be forced.
Similarly, if $R$ contains $-\vec{e}$,
then $-\vec{e}$ must be forced,
i.e., $\vec{e}$ must be forbidden.

For the other direction,
suppose $\vec{e}$ has endpoints
that are in different accessibility classes.
Let $G'$ be the graph obtained from $G$
by contracting all edges whose endpoints
lie in the same accessibility class,
let $N'$ be the number of vertices of $G'$,
and let $R_0$ be an orientation of $G'$
obtained by taking a $c$-orientation of $G$
and contracting all the directed edges
whose endpoints lie in the same accessibility class.
Note that $R_0$ is acyclic.
We will iteratively construct orientations $R_1,R_2,...$ of $G'$,
where each orientation $R_i$
is obtained from the preceding orientation $R_{i-1}$
by pushing down at a vertex.
Imagine the vertex set of $G'$ as being divided at each stage
into two classes, called ``red'' and ``blue,''
with the red class initially empty.
We will always be pushing down at a blue vertex of $G'$
and re-coloring it red.
More specifically, at the $i$th stage ($i \geq 1$)
we proceed as follows:
Pick an arbitrary blue vertex and proceed from there,
following edges of $R_{i-1}$.
Note that it is impossible
to go from a blue vertex to a red vertex,
since the directed edges of $R_{i-1}$
that connect red and blue vertices
always point away from the red (pushed-down) vertices.
Hence the path will only visit blue vertices.
Since $R_{i-1}$ was obtained from $R_0$ by repeatedly pushing down,
and since $R_0$ was acyclic, $R_{i-1}$ is acyclic as well;
hence the blue path must eventually get stuck at some blue vertex.
This vertex will be maximum, so that it is possible to push it down,
obtaining a new orientation $R_i$. Repeating this operation $N'$ times,
we eventually reach the orientation $R_{N'} = R_0$.
In the process of going from $R_0$ to $R_{N'}$,
every vertex of $G'$ has gotten pushed down;
hence every directed edge of $G'$ has gotten reversed twice.
In particular, we see that every edge of $G'$
can occur with both orientations.
Uncontracting the vertices of $G'$ to bring ourselves back to $G$,
we see that every edge of $G$ that joins different accessibility classes
occurs with both orientations in $\cR$.
\end{proof}

Before defining a partial ordering on $\cR$
(which will depend critically on our choice of $A^*$),
we need two extra ingredients:
one is a choice of a particular vertex $v^*$ in $A^*$
and the other is a real-valued function on the set of directed edges of $G$
that has certain properties.
To motivate these properties,
note that every $c$-orientation $R$ of $G$
is associated with an indicator function $F_R: \vec{E} \rightarrow \R$,
where
$$
F_R(\vec{e}) = \left\{ \begin{array}{cl}
        +1 & \mbox{if $\vec{e} \in R$}, \\
        0 & \mbox{if $\vec{e} \not\in R$ (i.e., $-\vec{e} \in R$)}.
        \end{array} \right.
$$
This function $F=F_R$ has three notable properties:
$0 \leq F(\vec{e}) \leq 1$ for all $\vec{e} \in \vec{E}$;
$F(\vec{e})+F(-\vec{e})=1$ for all $\vec{e} \in \vec{E}$;
and
$\sum_{\vec{e} \in C} F(\vec{e}) = \frac12(|C|+c(C))$
for every directed cycle $C$ in $G$.
(Actually the second property is a special case of the third.)

We would like a particular function $\overline{F}$
that satisfies these three properties
and some extra conditions as well:
$\overline{F}(\vec{e})$ should be $1$ when $\vec{e}$ is forced,
$0$ when $\vec{e}$ is forbidden,
and strictly between 0 and 1
when $\vec{e}$ is neither forced nor forbidden.
A simple way to find such an $\overline{F}$
is to average the $F_R$'s over all the (finitely many)
$c$-orientations $R$.
That is, we may let $\overline{F}(\vec{e})$ be
the probability that a $c$-orientation
chosen uniformly at random from $\cR$ contains $\vec{e}$.
In applications to specific graphs $G$ and circulations $c$,
there is typically a more natural function $\overline{F}$
that satisfies our extra property,
but for the purposes of proving general theorems
this choice of $\overline{F}$ suffices.
In any case, neither the choice of the vertex $v^* \in A^*$
nor the choice of the function $\overline{F}$
affects the partial ordering of $\cR$.

\begin{definition} \label{hf}
The {\bf height-function} $H_R$
of a $c$-orientation $R$
is the unique real-valued function $H$ on $V$
such that

(i) $H(v^*) = 0$ and

(ii) for $\vec{e} = (e,v,w)$, $H(w) - H(v)$ equals
$1 - \overline{F} (\vec{e})$ if $\vec{e} \in R$
and equals $- \overline{F} (\vec{e})$ otherwise.
\end{definition}

\begin{proposition} \label{PA}
The preceding definition is consistent and uniquely specifies $H_R$.
Moreover, every $H:V \rightarrow \R$ that satisfies (i) and (ii)
is equal to $H_R$ for some $R \in \cR$.
\end{proposition}

\begin{proof}
Since $G$ is connected,
we can find a directed path $\vec{e}_1,\vec{e}_2,...,\vec{e}_n$
from $v^*$ to any other vertex $v$ of $G$.
Then (i) and (ii) imply that
$H_R(w) = k - \overline{F}(\vec{e}_1)
- \overline{F}(\vec{e}_2)
- ...
- \overline{F}(\vec{e}_n)$,
where $k$ is the number of $i$ ($1 \leq i \leq n$)
such that $\vec{e}_i \in \cR$.
Hence the function $H_R$ is uniquely specified;
we need to check
that the definition is self-consistent.
Letting $C = \vec{e}_1,\vec{e}_2,...,\vec{e}_n$
be a directed cycle in $G$
we find that (ii) gives us
$\sum_{i=1}^{n} (H(v_i)-H(v_{i-1})) =
k - \sum_{i=1}^{n} \overline{F}(\vec{e}_i)$,
where $\vec{e}_k = (e_k,v_{k-1},v_k)$
and where $k$ is the number of $i$
with $\vec{e}_i \in \cR$;
we must check that this vanishes
(since $v_0=v_n$).
But $k$ equals the number of forward edges in the directed cycle $C$,
which is also equal to
$\sum_{i=1}^{n} \overline{F}(\vec{e}_i)$.
This proves the first sentence of Proposition \ref{PA}.

To construct an orientation of $G$ from a function $H$
satisfying (i) and (ii),
let $R$ be the set consisting of all forced directed edges
together with all directed edges $\vec{e}=(e,v,w)$
for which $H(w) > H(v)$.
$R$ is an orientation of $G$;
it is easily seen that
if $H=H_{R'}$ for some orientation $R'$,
we must have $R'=R$.
\end{proof}

Note that the condition (ii) in Definition \ref{hf} can be written
as $H_R(w)-H_R(v)=F_R(\vec{e})-{\overline F}(\vec{e})$,
where $F_R$ is the indicator function of $R$.

\begin{proposition} \label{PB}
For all $v,w$ adjacent in $G$,
$H_R(v)=H_R(w)$ if and only if
$v$ and $w$ belong to
the same accessibility class relative to $c$.
\end{proposition}

\begin{proof}
Let $e$ be an edge with endpoints $v$ and $w$.
If $v$ and $w$ are mutually accessible,
then the directed edge $\vec{e} = (e,v,w)$
is either forced or forbidden
so that $F_R(\vec{e}) = {\overline F}(\vec{e})$
and $H_R(w) = H_R(v)$.
If $v$ and $w$ are not mutually accessible,
$\vec{e}$ is neither forced nor forbidden,
so that either
$H_R(w)-H_R(v) = 1 - \overline{F}(\vec{e}) > 0$
or $H_R(w)-H_R(v) = - \overline{F}(\vec{e}) < 0$.
\end{proof}

\begin{proposition}
\label{PC}
For all $R,R' \in \cR$ and for all $v \in V$,
$H_R(v) - H_{R'}(v)$ is an integer.
\end{proposition}

\begin{proof}
It trivially holds for $v=v^*$ (by stipulation (i)
in the definition of $H_R$),
and since the two alternatives on the right hand side of (ii)
differ by 1, it holds by induction for all $v$.
\end{proof}

\begin{definition}
For each $v \in V$ let $m(v) = \min_{R \in \cR} H_R(v)$.
\end{definition}

\begin{proposition}
\label{PCa}
Suppose $v$ and $w$ are adjacent vertices
in different accessibility classes of $G$,
with $m(v) < m(w)$. Then for $R \in \cR$,
the only possible values for the pair $(H_R(v),H_R(w))$
are of the form $(m(v)+i,m(w)+j)$ with $j=i$ or $j=i-1$.
\end{proposition}

\begin{proof}
Note that the assumption $m(v) < m(w)$
entails no loss of generality,
since $m(v) \neq m(w)$
whenever $v$ and $w$ are in different accessibility classes
(Proposition \ref{PB}).
I claim that $m(w)$ lies strictly between $m(v)$ and $m(v)+1$.
Indeed, our hypothesis $m(v)<m(w)$ gives us one bound.
To prove the other,
let $R$ be an orientation for which $H_R(v)=m(v)$.
Then $H_R(w)$ is either $m(v)-{\overline F}(\vec{e})$
or $m(v)+1-{\overline F}(\vec{e})$,
where $\vec{e}$ is the directed edge from $v$ to $w$.
Since $\vec{e}$ is neither forced nor forbidden,
$m(v)-{\overline F}(\vec{e}) < m(v)$
and $m(v)+1-{\overline F}(\vec{e}) < m(v)+1$.
Since $H_R(w) \geq m(w) > m(v)$,
we cannot have $H_R(w) = m(v) - F(\vec{e})$;
hence $H_R(w) = m(v)+1-{\overline F}(\vec{e})$.
Therefore in fact we have
$m(w)=H_R(w)=m(v)+1-{\overline F}(\vec{e})<m(v)+1$.

Now let $R$ be any orientation of $G$.
By Proposition \ref{PC},
$H_R(v)=m(v)+i$ and $H_R(w)=m(w)+j$ for some $i,j$.
We must either have $H_R(w)-H_R(v)=1-{\overline F}(\vec{e})=m(w)-m(v)$
or $H_R(w)-H_R(v)=-{\overline F}(\vec{e})=m(w)-m(v)-1$.
In the first case, $j=i$; in the second case, $j=i-1$.
\end{proof}

\begin{proposition} \label{PD}
If $H_1$ and $H_2$ satisfy (i) and (ii),
then so do their {\bf meet} $H_1 \meet H_2$
and their {\bf join} $H_1 \join H_2$, where
$$
(H_1 \meet H_2) (v) = \min(H_1(v),H_2(v)) ,
$$
$$
(H_1 \join H_2) (v) = \max(H_1(v),H_2(v)) .
$$
\end{proposition}

\begin{proof}
It is clear that (i) is satisfied in each case.
To prove (ii),
note that if $v$ and $w$ are mutually accessible,
then the possible values for $(H_R(v),H_R(w))$
as $R$ varies in $\cR$ are
$$(r,r),(r+1,r+1),(r+2,r+2),...$$
where $r=m(v)$,
whereas if $v$ and $w$ are not mutually accessible,
then the possible values for $(H_R(v),H_R(w))$ are
$$(r,s),(r+1,s),(r+1,s+1),(r+2,s+1),...$$
where $r=m(v)$ and $s=m(v)$ (Proposition \ref{PCa}).
Since in either situation the allowed pairs are linearly ordered componentwise,
any ordered pair of numbers that is of the form
$(H_1(v),H_1(w)) \join (H_2(v),H_2(w))$
or
$(H_1(v),$ $H_1(w)) \meet (H_2(v),H_2(w))$
is also a pair on the list, and hence satisfies (ii).
\end{proof}

\begin{definition}
If $R$ and $S$ are $c$-orientations, say $R \succeq S$
if and only if $H_R(v) \geq H_S(v)$ for all $v \in V$.
\end{definition}

Clearly, $(\cR, \meet, \join, \succeq)$ is a distributive lattice.
We define $\succ$, $\preceq$ and $\prec$
in terms of $\succeq$ in the usual way.
Note that $R \succeq S$ if and only if
for every vertex $v$ of $G$,
and for every path in the graph from $v^*$ to $v$,
the number of forward edges along the path
is at least as great for $R$ as for $S$.

We must now show that $R$ covers $S$
in the partial order $(\cR, \succeq)$
if and only if $S$ is obtained from $R$
by a single application of pushing-down.

\begin{proposition} \label{PE}
Let $R$ be a $c$-orientation,
and let $A$ be an accessibility class of $G$
(relative to $c$).
Then $A$ is maximal under $R$
(that is, all directed edges between $A$ and $A^c$
point towards $A$)
if and only if
$A$ is a ``mesa'' for the height-function $H_R$;
that is, $H_R$ is constant on $A$,
and has strictly smaller value for every vertex in $A^c$
adjacent to a vertex in $A$.
\end{proposition}

\begin{proof}
By Proposition \ref{PB} we know
that $H_R$ is constant on accessibility classes;
say $H_R(v)=\alpha$ for all $v \in A$.
For each $w \not\in A$ adjacent to $A$,
we have either $H_R(w) < \alpha$
or $H_R(w) > \alpha$,
according to whether the edge(s) connecting $w$ to $A$
point towards $A$ or away from $A$.
(Note that it is impossible for there to be
edges joining $w$ to $A$ with both orientations,
since that would force $w \in A$.)
Thus $A$ is maximal if and only if
$H_R(w) > \alpha$ for all $w \not \in A$ adjacent to $A$.
\end{proof}

\begin{proposition} \label{PF}
If $A$ is a maximal accessibility class in $\cR$,
and $S$ is obtained from $R$ by pushing down $A$, then
$H_S(v)$ equals $H_R(v)-1$ if $v \in A$
and equals $H_R(v)$ otherwise.
\end{proposition}

\begin{proof}
Easy.
\end{proof}

\begin{proposition} \label{PG}
Let $R,S$ be $c$-orientations with $R \not\preceq S$,
so that $H_R(v) > H_S(v)$ for some $v$, and let
$d = \max \{ H_R(w)-H_S(w): w \in V \} > 0$.
Let
$D = \{ w \in V: H_R(w) - H_S(w) = d \}$,
which by Proposition \ref{PB} is a union of accessibility classes,
and let $\hat{A}$ be an accessibility class in $D$
for which $H_R(\hat{A})$ is as great as possible.
Then $\hat{A}$ is maximal in $R$.
Moreover, if $R \succ S$,
then the orientation $R'$ that results from $R$
by pushing down $\hat{A}$ satisfies $R' \succeq S$.
\end{proposition}

\begin{proof}
Since $H_R(v)-H_S(v)=0$ for all $v \in A^*$,
$D$ is disjoint from $A^*$.
Let $w$ be a vertex in $\hat{A}^c$ adjacent to $\hat{A}$.
If $w \in D$, then $H_R(\hat{A}) > H_R(w)$,
since $H_R$ is maximized within $D$ on the set $\hat{A}$
($H_R(w)=H_R(\hat{A})$ is impossible since $w \not\in \hat{A}$).
On the other hand,
if $w$ does not belong to $D$,
so that $H_R(w)-H_S(w) < d$,
then $H_R(w)-H_S(w)=d-1$,
and our earlier observation about the possible values
of the pair $(H(v),H(w))$ (for $H$ any height-function)
gives us $H_R(\hat{A}) > H_R(w) \geq H_S(w) > H_S(\hat{A})$.
In either case, $H_R(\hat{A}) > H_R(w)$,
and since this holds for all $w$ adjacent to $\hat{A}$,
all the arrows between $\hat{A}$ and its complement point towards $\hat{A}$,
implying that $\hat{A}$ is maximal in $R$.
To prove the second claim,
note that the inequality $H_R(\hat{A}) > H_S(\hat{A})$
implies $H_R(\hat{A}) \geq H_S(\hat{A})+1$.
Hence $H_{R'}(w) = H_R(w) \geq H_S(w)$ for all $w \not\in \hat{A}$ and
$H_{R'}(\hat{A}) = H_R(\hat{A})-1 \geq H_S(\hat{A})$.
\end{proof}

\begin{proposition} \label{PH}
$R$ covers $S$ if and only if
$S$ is obtainable from $R$
by a single pushing-down operation
on an accessibility class $A$ disjoint from $A^*$.
\end{proposition}

\begin{proof}
The backward direction follows immediately from Proposition \ref{PF},
since heights at a fixed vertex $v$
can only vary by integers as $R$ varies
and since $H_R$ must be constant on accessibility classes.
In the forward direction,
suppose $R$ covers $S$ in $(\cR, \succeq)$.
By Proposition \ref{PG},
there exists $R' \succeq S$ obtained from $R$
by pushing down.
Since $R \succ R'$,
we must have $R'=S$.
\end{proof}

This interpretation of pushing down gives us the following result:

\begin{proposition} \label{PI}
There is a unique $c$-orientation of $G$
which has no maximum away from $A^*$.
\end{proposition}

\begin{proof}
Any $c$-orientation of $G$ with no maximum away from $A^*$
is a $c$-orientation that is minimal in $\cR$,
and vice versa.
Since $\cR$ is a lattice,
the orientation is unique.
Note that the orientation must necessarily have a maximum at $A^*$,
since its height-function must achieve its maximum someplace.
\end{proof}

Call this orientation $R_{\hat{0}}$.
Similarly, we define $R_{\hat{1}}$
as the unique $c$-orientation of $G$
which has no minimum away from $A^*$;
$R_{\hat{1}}$ is the maximum element of $\cR$.

We can now prove Theorem \ref{thm1}.

\begin{proof}
Summarizing our results, we find that
$(\cR, \meet, \join, \succeq)$
is a distributive lattice
in which an element $R$ covers an element $S$
if and only if $S$ is obtained by pushing down $R$
at a maximal accessibility class other than $A^*$,
with minimum element $R_{\hat{0}}$ and maximum element $R_{\hat{1}}$.
This completes the proof of Theorem \ref{thm1}.
\end{proof}

We now turn to examples.

\begin{example} \label{E11}
Let $G$ be a cycle of length $n$,
and let $\cR$ be the set of orientations of $G$
in which exactly $k$ of the edges of $G$ are directed counterclockwise.
Here we put $\overline{F}(\vec{e}) = \frac{k}{n}$
for every counterclockwise edge and $\frac{n-k}{n}$ for every clockwise edge.
(See Figure 2 for the case $n=4$, $k=2$.)
In general, $\cR$ is isomorphic to the finite Young lattice ${\rm L}(k,n-k)$.
To see this, number the edges of $G$ from 1 to $n$ in clockwise order
(beginning and ending with $v^*$),
and create a lattice path from $(k,0)$ to $(0,n-k)$
whose $i$th edge goes upward or leftward
according to whether the $i$th edge of $P$
runs clockwise or counterclockwise.
The set of grid squares in the first quadrant
to the left of and below the lattice path is a Young diagram,
and pushing down amounts to removing a square.
Thus, the $c$-orientations of $G$ correspond to number-partitions
having at most $k$ parts, each of size at most $n-k$.
(For a general treatment of number-partitions, see \cite{Andr}.)
\end{example}

\begin{example} \label{E12}
Let $G$ be a path composed of $n$ horizontal edges,
labeled $e_1,e_2,...,e_n$ from left to right,
with $v^*$ the rightmost vertex.
Since $G$ has no cycles,
we let $\cR$ be the set of all orientations of $G$.
For any orientation $R$ of $G$,
let $I(R) \subseteq \{1,2,...,n\}$
be the set of $i$ ($1 \leq i \leq n$)
such that $e_i$ is oriented from right to left in $R$.
Then $R$ covers $S$ if and only if $I(S)$ is obtained from $I(R)$
by reducing one of its elements by 1
(with the stipulation that this reduction is only permitted
if the reduced value was not already in $R$,
and with the special arrangement that
the number 1 is ``reduced'' by being simply deleted).
Thus $\cR$ is isomorphic to the inclusion-ordering on number-partitions
with distinct parts of size at most $n$.
(More generally, if $G$ is any tree
with $n$ edges and $v^*$ is any vertex of $G$,
$\cR$ will be a poset with $2^n$ elements.)
\end{example}

For our next example,
we need to broaden the context of our construction somewhat.
Let $\partial G$ denote any non-empty connected subgraph of $G$,
to be called the boundary of $G$.
Let $\partial V$ denote the set of vertices of $\partial G$;
assume $v^* \in \partial V$.
If we let $\cR$ be the set of all $c$-orientations of $G$
that have prescribed orientations on all the edges of $\partial G$,
then the proof of the preceding propositions can be adapted
to show that every orientation in $\cR$
can be obtained from every other
by means of pushing-down operations
in which the accessibility classes of $G$ that get pushed down
are disjoint from $\partial V$.

\begin{example} \label{E13}
Let $G$ be the square grid with $(n+1)^2$ vertices ($n \geq 1$),
let $\partial G$ be the boundary of $G$ in the usual sense
(a subgraph with $4n$ vertices and $4n$ edges),
let $c$ be the circulation that assigns circulation zero
to every cycle, and let $\cR$ be the set of $c$-orientations of $G$
in which all edges on the top/bottom/left/right side of the boundary
point rightward/leftward/
downward/upward.
One such $R \in \cR$ is shown in Figure 6(a) for $n=3$.

\begin{figure}
\begin{center}
\includegraphics[width=6.0in]{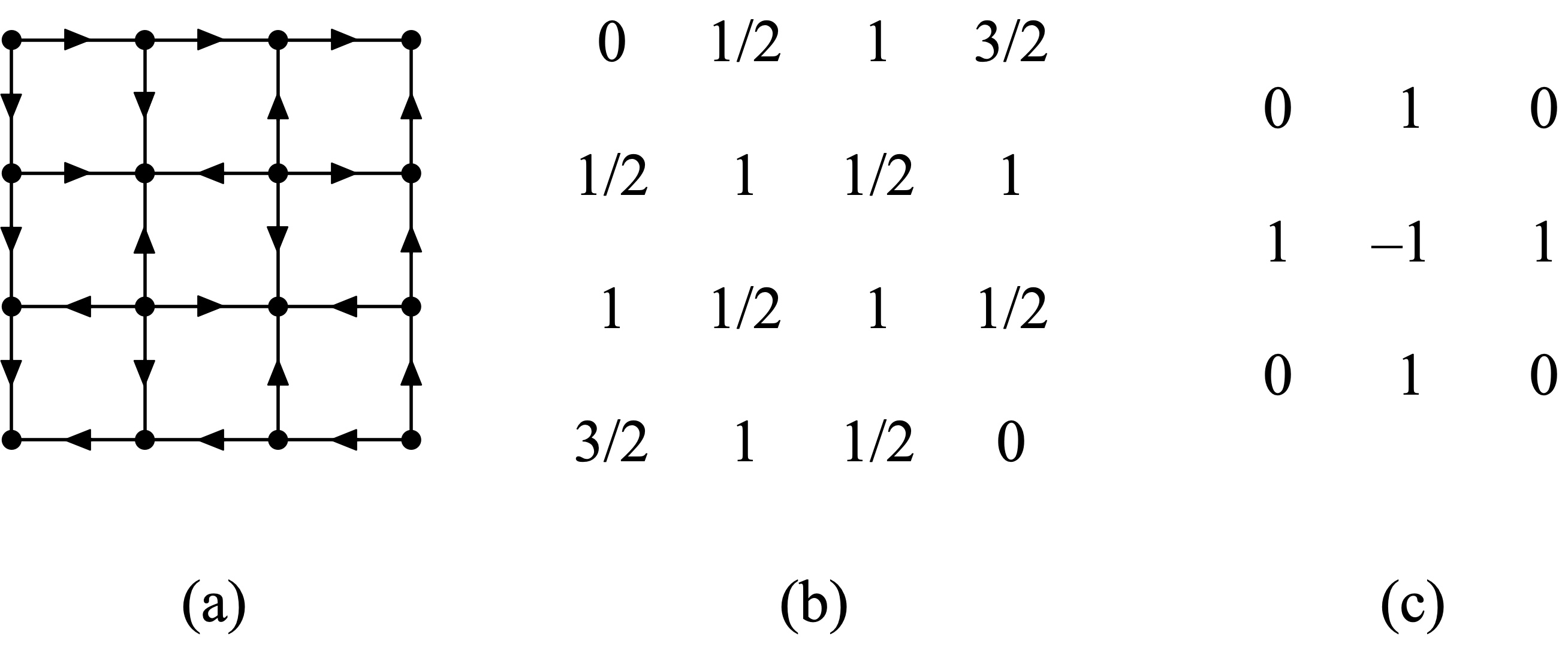}
\end{center}
\caption{Alternating sign matrices.}
\end{figure}

If we take $\overline{F}(\vec{e}) = \frac12$ for all $\vec{e}$,
with $v^*$ the vertex at the upper left,
we get the height-function $H_R$ shown in Figure 6(b).
Finally, if for each grid-square in $G$
we record the signed sum $-h_1+h_2+h_3-h_4$
where $h_1,h_2,h_3,h_4$ are
the heights $H_R$ of the upper left, upper right,
lower left, and lower right corners of the grid-square, respectively,
we get an $n$-by-$n$ array of $+1$'s, $-1$'s, and $0$'s, as shown in Figure 6(c).
This array has the property that in each row and column,
the nonzero entries alternate in sign, beginning and ending with $+1$.
Such a matrix is called an {\bf alternating-sign matrix}
(defined by Mills, Robbins, and Rumsey~\cite{MRR1};
for a very readable overview, see~\cite{Robb}).
It is not hard to show that every $n$-by-$n$ alternating-sign matrix
corresponds to a unique orientation of $G$ in $\cR$
(see Elkies et al.~\cite{EKLP} for details).
Our lattice structure on $\cR$ is
a natural partial ordering of the set of alternating-sign matrices,
and indeed has been used numerous times in the literature.
One indication of its naturalness is the observation that
if we restrict this ordering to those alternating-sign matrices
in which no $-1$'s occur (the permutation matrices),
then we obtain the weak Bruhat order on permutations.
Conversely, our ordering of alternating sign matrices
can be construed in order-theoretic terms as the
MacNeille completion of the weak Bruhat order; see~\cite{BS}.
\end{example}

For the rest of this section,
we will assume that $c$ is acyclic,
and that choosing an orientation for the edges in $\partial G$
does not force any of the orientations of the other edges.

Since $(\cR, \succeq)$ is a distributive lattice,
it can be viewed as the lattice of order-ideals of a poset $P$,
ordered by inclusion.
$P$ can be abstractly described as
the poset of join-irreducibles of $\cR$
under the inherited partial ordering,
but here is a more concrete picture of $P$.
For each vertex $v \not\in \partial V$,
let $m(v) = \min_R H_R(v)$,
$M(v) = \max_R H_R(v)$,
and $D(v)=M(v)-m(v) \geq 1$,
so that the allowed values of $H_R(v)$
are $m(v),m(v)+1,...,m(v)+D(v)=M(v)$.
To each $v \not\in \partial V$
and each $1 \leq i \leq D(v)$,
there corresponds an element of $P$,
which we can denote by $(v,i)$;
every element of $P$ is of this form.
$(v,i)$ covers $(w,j)$ in $P$
when $v$ and $w$ are adjacent in $G$
and $0 < (m(v)+i) - (m(w)+j) < 1$.

A different geometric picture of $(\cR, \succeq)$
arises from looking at the
\boldmath $c${\bf -orientation polytope}
\unboldmath
associated with the circulation $c$.
As discussed earlier, each orientation $R$ of $G$
can be associated with an indicator function $F_R$.
We can view the set of real-valued functions
from $\vec{E}$ to $\R$ as a Euclidean space
and consider the functions $F_R$ as points in that space.
We define the $c$-orientation polytope as the convex hull of these points.
Note that it actually lies in the subspace of $\R^{\vec{E}}$
consisting of those $F$
that satisfy $\sum_{\vec{e} \in C} F(\vec{e}) = \frac12(|C|+c(C))$
for all directed cycles $C$.

\begin{proposition} \label{PJ}
If $S$ is obtained from $R$ by pushing down a vertex $v$,
then there is an edge of the $c$-orientation polytope
between $R$ and $S$.
\end{proposition}

\begin{proof}
First note that every $F_R$ is a vector containing $|E|$ 0's
and $|E|$ 1's;
hence none is a convex combination of the others
and each is on the boundary of the $c$-orientation polytope.
To see that there is an edge from $F_R$ to $F_S$,
consider the (affine) subspace
consisting of all $F \in \R^{\vec{E}}$
that agree with $F_R$ and $F_S$
on all directed edges common to $R$ and $S$,
and intersect it with the subspace of $F$'s
satisfying the identity
$\sum_{\vec{e} \in C} F(\vec{e}) = \frac12(|C|+c(C))$.
We get a line and nothing more,
since the values of $F$ on all directed edges not involving $v$
are fixed on the first subspace
and the differences $F(\vec{e}_i) - F(\vec{e}_j)$
are fixed on the intersection with the second subspace
(where $\vec{e}_i$ and $\vec{e}_j$
are any two directed edges with initial vertex $v$).
\end{proof}

\begin{proposition} \label{PK}
There exists an affine function $\Phi$ on $\R^{\vec{E}}$
such that $\Phi(F_R)$ is the rank of $R$ in the poset $(\cR, \succeq)$.
\end{proposition}

\begin{proof}
We seek a function of the form
$$\Phi(F) = \sum_{\vec{e} \in \vec{E}} a_{\vec{e}} F(\vec{e}) + a_0$$
for coefficients $a_{\vec{e}}$ and a constant term $a_0$.
It suffices to find $a_{\vec{e}}$'s so that
every pushing-down operation
decreases $\sum_{\vec{e}} a_{\vec{e}} F(\vec{e})$ by 1.
That is, if $v \neq v^*$ is a vertex of $G$
and $\vec{E}(v)$ is the set of directed edges of $G$
with initial vertex $v$,
then we want
$$\sum_{\vec{e} \in \vec{E}(v)} a_{-\vec{e}} =
\left( \sum_{\vec{e} \in \vec{E}(v)} a_{\vec{e}} \right) + 1 \ .$$
Write this as
$$\sum_{\vec{e} \in \vec{E}(v)}
(a_{-\vec{e}} - a_{\vec{e}}) = 1 \ .$$
To see if there are any degeneracies in this system of equations,
consider all linear combinations of the form
$$\sum_{v \in V} \left(
b_v \sum_{\vec{e} \in \vec{E}(v)}
(a_{-\vec{e}} - a_{\vec{e}}) \right)$$
with $b_{v^*}$ required to be zero.
To see when this expression vanishes identically
(that is, to see when the net coefficient of each $a_{\vec{e}}$ is zero),
rewrite the expression as
$$\sum_{\vec{e}=(e,v,w) \in \vec{E}} (b_w - b_v) a_{\vec{e}} \ . $$
This vanishes identically only when $b_w = b_v$
for all $v,w$ adjacent in $G$;
since $G$ is connected,
this means that $b_v = b_{v^*} = 0$ for all $v$.
Hence our inhomogeneous system of equations contains no inconsistencies
and has at least one solution.
\end{proof}

One curious feature of Theorem \ref{thm1}
is that it gives us $N$ different ways
to view $\cR$ as a distributive lattice,
where $N$ is the number of vertices of $G$.
Let us write $\cR_v$ for the partial ordering
on the $c$-orientations of $G$
defined using $v$ as the special vertex.
If we let $P_v$ be the poset of join-irreducibles of $\cR_v$,
then we obtain commuting bijections
$J(P_v) \rightarrow J(P_w)$
for all $v,w$,
where $J(P)$ denotes the lattice of order ideals of $P$.
These are non-trivial bijections,
in that the underlying $P_v$'s can be non-isomorphic.

Abstractly, what is going on is best understood
in the context of the infinite distributive lattice $\hat{\cR}$
consisting of all functions $H:V \rightarrow \R$ for which
$H(v^*) \in \Z$ and $H(w)-H(v)$ equals either
$1-\overline{F}(\vec{e})$ or
$-\overline{F}(\vec{e})$
for all directed edges $\vec{e}=(e,v,w) \in \vec{E}$.
Note that $v^*$ no longer plays a truly special role;
the structure of this lattice is independent of our choice of $v^*$.
For instance, in the case where $G$ is a path with two edges,
the infinite lattice
is as shown in Figure 7(a);
its join-irreducibles are ordered in the fashion
shown in Figure 7(b).
The lattice admits a translational symmetry,
which corresponds to adding 1 to every height.
When we pass from $\cR$ to some $\cR_v$,
we are in a sense modding out by this translation symmetry.
Here is the relevant general fact:

\begin{figure}
\begin{center}
\includegraphics[width=3in]{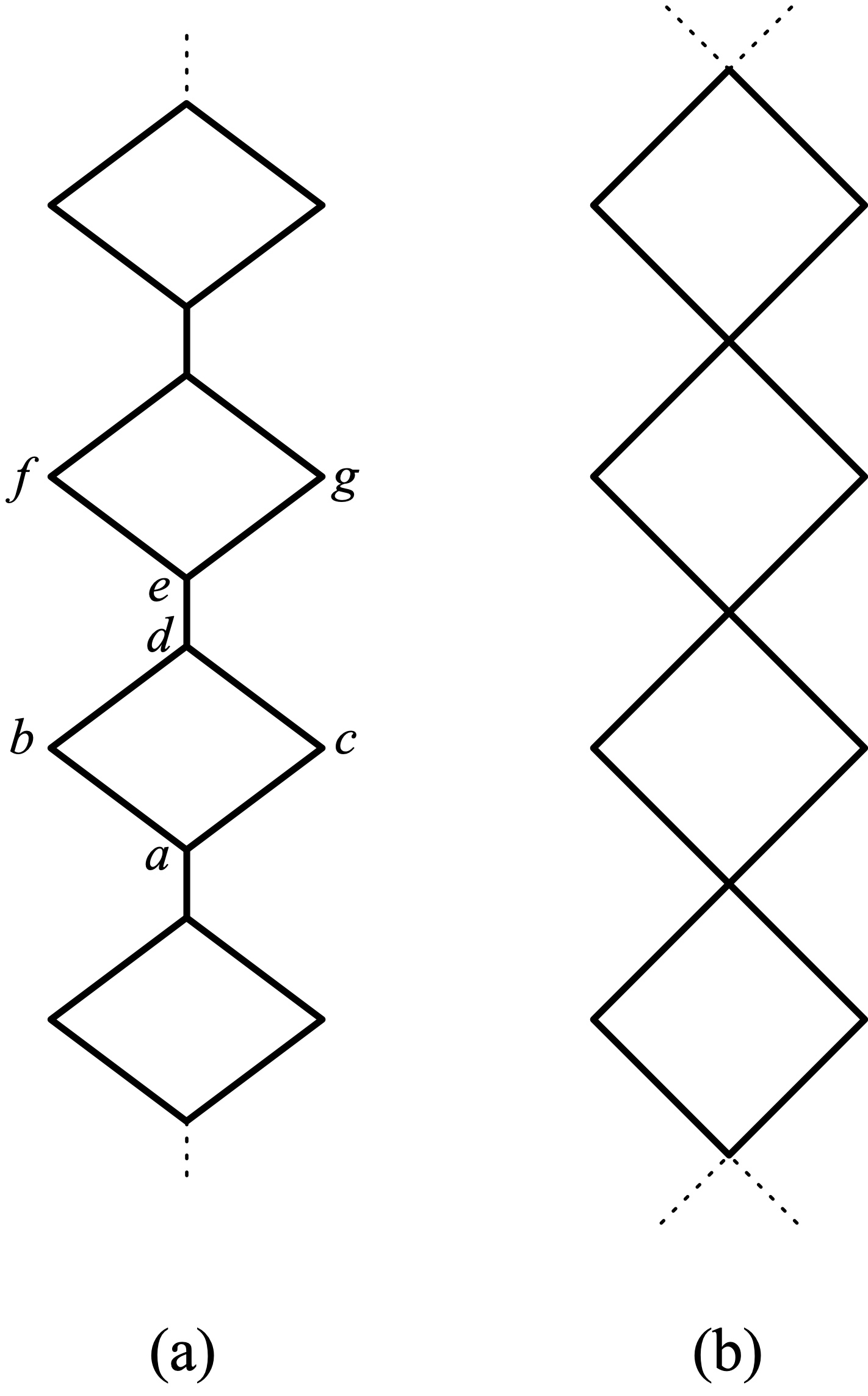}
\end{center}
\caption{An infinite lattice and its join-irreducibles.}
\end{figure}

\begin{proposition} \label{PL}
Let $\hat{L}$ be an infinite distributive lattice
with a free action of $\Z$,
such that there are finitely many orbits under the action,
each of which is a doubly-infinite chain.
Let $\hat{P}$ be the poset of join-irreducibles
of $\hat{L}$,
which carries an induced action of $\Z$
under which there are finitely many orbits,
each of which is a doubly-infinite chain.
Define an {\bf interval} of the lattice
as a set of the form $[a,b] = \{x: a \preceq x \preceq b\}$.
If $I$ is any finite interval in $\hat{L}$
whose translates partition $\hat{L}$,
then $I$ is a finite distributive lattice
and there is a bijection between $I$ and
the set of orbits of $L$ under $\Z$.
\end{proposition}

\begin{proof}
Easy.
\end{proof}

Call such an interval $I$
a {\bf finite quotient} of the infinite lattice $\hat{L}$.
Then, returning to our original context
($c$-orientations of graphs),
we see that
the different $P_v$'s associated with different vertices $v$
all arise as the poset of join-irreducibles
of different finite quotients of $\hat{\cR}$.
For instance, in Figure 7(a),
the intervals $\{a,b,c,d\}$, $\{b,d,e,g\}$, and $\{c,d,e,f\}$
are all quotients of the lattice $\hat{L}$.
The first is isomorphic to the lattice of order-ideals
of a two-element anti-chain;
the second and third are isomorphic to
the lattice of order-ideals of a three-element chain.

It would be interesting to know when two abstract posets $P_1,P_2$
can arise from two finite quotients of the same $\hat{L}$.
For instance,
if $P$ is the product of a two-element chain with itself,
can the construction of Proposition \ref{PL} be used
to give a bijection between $J(P)$ and $J(P')$,
where $P'$ is some poset not isomorphic to $P$?

As promised earlier,
I will now apply the theory of height-functions
to the question of bounding the number of pushing-down operations
required to get from one $c$-orientation to another.

\begin{proposition} \label{PM}
It is possible to get from any $c$-orientation of a graph $G$
to any other in at most $N(N-1)/2$ push-down operations,
where $N=|V|$.
\end{proposition}

\begin{proof}
Let $R,S$ be two $c$-orientations of $G$.
Take a vertex $v^*$ at random,
and calculate the height-functions $H_R$ and $H_S$
in the ordering $\cR_{v^*}$.
Let $v'$ be the vertex $v$ of $G$
for which $H_R(v)-H_S(v)$
is as small (or as negative) as possible.
If we now switch over to using $v'$ as our special vertex
rather than $v^*$,
we are effectively sliding the two height-functions
in the vertical direction
until they touch (at $v'$) but do not cross.
In the partial ordering $\cR_{v'}$,
we have $R \succeq S$.
Hence, by repeated application of Proposition \ref{PG},
one sees that it is possible to convert $R$ into $S$
in $\sum_v |H_R(v)-H_S(v)|$ steps,
each of which brings $H(v)$ closer to $H_S(v)$
for one vertex $v$
while leaving the rest of the height-function alone.
Therefore the number of moves required is at most
$\sum |H_{R_{\hat{1}}}(v) - H_{R_{\hat{0}}}(v)|
= \sum D(v)$,
where $D(v) = \max_H H(v) - \min_H H(v)$ as before.
Let $k = \max_v D(v)$.
Note that for $v$ and $w$ adjacent,
$D(v)$ and $D(w)$ differ by at most 1,
and that $D(v')=0$.
Hence there must exist $v_i \in V$ with $D(v_i) = i$
for $i=1,2,...,k$.
$D(v) \leq k$ for all $v \not\in \{v', v_1, v_2, ..., v_k\}$,
so $\sum D(v) \leq 1 + 2 + ... + (k-1) + k + k + ... + k
= k(k-1)/2 + k(N-k) = k(2N-1-k)/2 \leq N(N-1)/2$.
\end{proof}

In the other direction, we can show that the bound $N(N-1)/2$ 
is best-possible:

\begin{proposition} \label{PN}
For every $N$,
there exists a graph $G$ on $N$ vertices
and orientations $R$ and $S$
with the same circulation around all directed cycles,
such that $S$ cannot be obtained from $R$
in fewer than $N(N-1)/2$ push-down operations.
\end{proposition}

\begin{proof}
Take $G$ as in Example \ref{E12}.
Let $R$ and $S$ be the orientations of $G$
in which all edges are directed
leftward (away from $v^*$) or rightward (towards $v^*$), respectively.
A pushing-down operation can be understood as
the operation of sliding a left-pointing edge toward the left
until it ``slides off the edge'',
or introducing a new left-pointing edge at $v^*$.
In order to convert $R$ into $S$,
all the left-pointing edges must be slid off of $G$,
so $1+2+...+(N-1)=N(N-1)/2$ operations are needed.
\end{proof}

\section{Matchings of Bipartite Graphs}

Suppose $G$ is a connected bipartite plane graph
with vertex-set $V$
and $d$ is a function from $V$ to the non-negative integers
such that $G$ has at least one $d$-factor.
Assume that every edge of $G$
belongs to some $d$-factors but not to others.
Fix an alternating coloring of the vertices of $G$
(black and white).
I will use duality for plane graphs
to apply the results of the preceding section
and get a family of partial orderings
on the $d$-factors of $G$.

For the reader's convenience, I restate Theorem \ref{thm2}.

\bigskip

\noindent
{\bf Theorem \ref{thm2}:} {\em Let $\cM$ be the (non-empty) set of $d$-factors
of a finite graph $G$ drawn on the sphere.
If we say that one $d$-factor $M$ covers another $d$-factor $N$
exactly when $N$ is obtained from $M$ by twisting down
at a face other than $f^*$,
then the covering relation makes $\cM$ into
a distributive lattice.}

\begin{proof}
Let $G^\perp$ be the dual of $G$
with vertex ${f^*}{^\perp}$ associated with face $f^*$.
The 2-coloring of the vertices of $G$
gives rise to a 2-coloring of the faces of $G^\perp$.
Define the {\bf standard orientation} of $G^\perp$
as the one in which all edges
are directed so as to encircle white faces of $G^\perp$
in the counterclockwise direction
and black faces of $G^\perp$ in the clockwise direction.
(See Figure 8.)

\begin{figure}
\begin{center}
\includegraphics[width=0.8\textwidth]{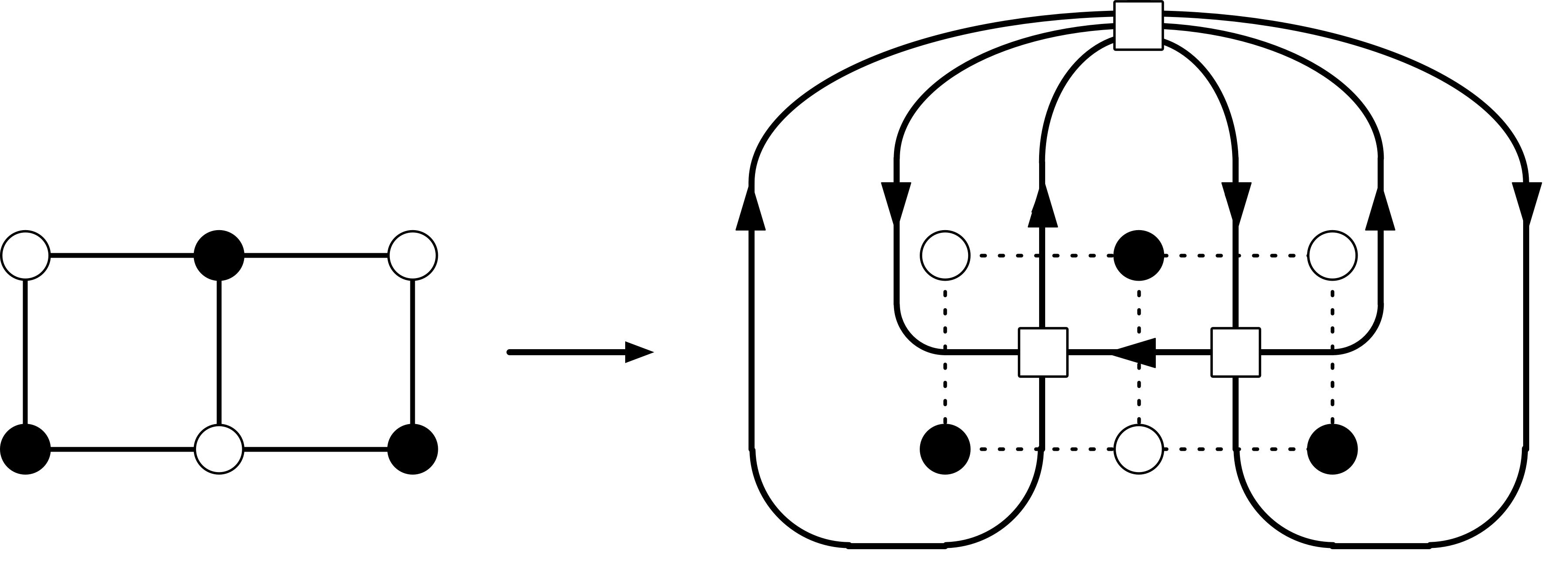}
\end{center}
\caption{The standard orientation of the dual.}
\end{figure}

Given a $d$-factor $M$,
we define an orientation $R_M$ of $G^\perp$
by giving an edge of $G^\perp$
its non-standard orientation
if the associated edge of $G$ belongs to $M$
and giving the edge its standard orientation otherwise.
If $C$ is an elementary cycle of $G^\perp$
that circles a white vertex of $G$ in the counterclockwise sense,
then the circulation of $R_M$ around $C$ is $\deg v - d(v)$;
if $C$ is an elementary cycle of $G^\perp$
that circles a black vertex of $G$ in the counterclockwise sense,
then the circulation of $R_M$ around $C$ is $-(\deg v - d(v))$;
and if $C$ is any other cycle of $G^\perp$,
the circulation of $R_M$ around $C$ is derivable
as a combination of circulations around
black and white vertices of $G$.
Thus every $R_M$ has the same circulation,
which we denote by $c$.
In the other direction, it is easy to check
that every $c$-orientation $R$ of $G^\perp$ singles out
the $d(v)$ edges at each vertex $v$ of $G$
that have their non-standard orientation
and thereby specifies a $d$-factor $M$ of $G$.
This gives a bijection between $\cR$ and $\cM$.
Moreover, positive/negative alternating cycles in $M$
correspond to maximal/minimal vertices in $R$, respectively,
twisting down a face in $M$
corresponds to pushing down a vertex in $R$,
and excluding the cycle bounding $f^*$
corresponds to excluding the vertex ${f^*}{^\perp}$.
Thus, Theorem \ref{thm1}, applied to $G^\perp$,
yields Theorem \ref{thm2}, applied to $G$.
\end{proof}

Note that under the bijection between $\cM$ and $\cR$,
one gets a height-function defined on the faces of $G$.
Since this height-function must take on
its maximum and minimum values somewhere,
we see that every $d$-factor of $G$
must have at least one positive alternating cycle
and at least one negative alternating cycle.

We can describe the lattice structure of $\cM$
without reference to the dual graph.
For each face $f \neq f^*$ of $G$,
choose a sequence of faces $f^*=f_0$, $f_1$, ..., $f_n=f$
such that every pair of consecutive faces $f_{i-1},f_i$
share an edge $e_i$
and such that, as one moves from $f_i$ to $f_{i+1}$,
the black vertex of $e_i$ is on the left.
Define the height of the face $f$
relative to a $d$-factor $M$
as the number of $e_i$'s that belong to $M$.
(As an alternative, one can use face-paths from $f$ to $f^*$
in which it is not required
that the black vertex of a shared edge be on the left,
but then the height of $f$ must be defined as
the number of even $e_i$'s in $M$
minus the number of odd $e_i$'s in $M$,
where an $e_i$ is called even or odd
according to whether the black vertex
is on the left or right.)
In any case, this height-function is the same as
the one we introduced for orientations of $G^\perp$,
modulo some renormalization at each face of $G$.
In particular, the lattice structure on $\cM$
given by Theorem \ref{thm2}
is the one induced by the lattice structure
on the set of these height-functions.

\begin{example} \label{E21}
Consider a hexagon that has all internal angles equal to 120 degrees,
with sides having integer lengths $a,b,c,a,b,c$, respectively.
We can imagine it as being divided into $2ab+2ac+2bc$
equilateral triangles of side 1,
half of them pointing upward
and half of them pointing downward.
Define a {\bf lozenge} as a rhombus of side 1
having angles of 60 and 120 degrees.
A tiling of the hexagon by lozenges
corresponds to a pairing between
the upward-pointing and downward-pointing triangles
in which paired triangles must have an edge in common.
That is to say, a tiling corresponds to a matching
of a certain bipartite graph $G$ with $2ab+2ac+2bc$ vertices.
(See Figure 9.)
Note that an alternating cycle of $G$ corresponds 
to a unit hexagon in the tiling composed of three lozenges, 
and that twisting such a cycle corresponds 
to rotating the hexagon by 180 degrees.

\begin{figure}
\begin{center}
\includegraphics[width=0.8\textwidth]{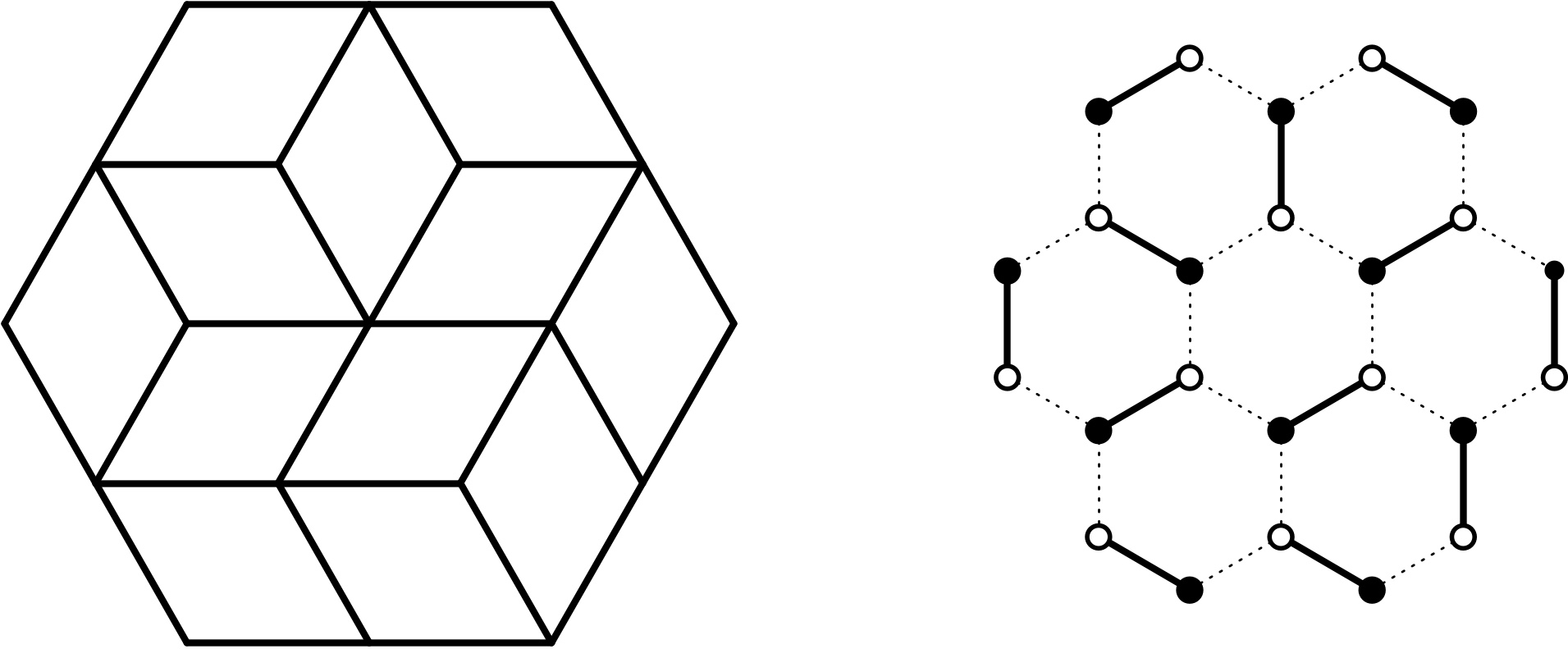}
\end{center}
\caption{Lozenge tilings and matchings.}
\end{figure}

If we color vertices black or white
according to whether the associated triangle
points upward or downward
and we let $f^*$ be the external face of $G$,
then we obtain a partial ordering on the set of matchings of $G$,
which corresponds to a partial ordering on
the set of tilings of the hexagon.
(Figure 10 shows the case $a=1$, $b=2$, $c=2$.)

\begin{figure}
\begin{center}
\includegraphics[width=0.5\textwidth]{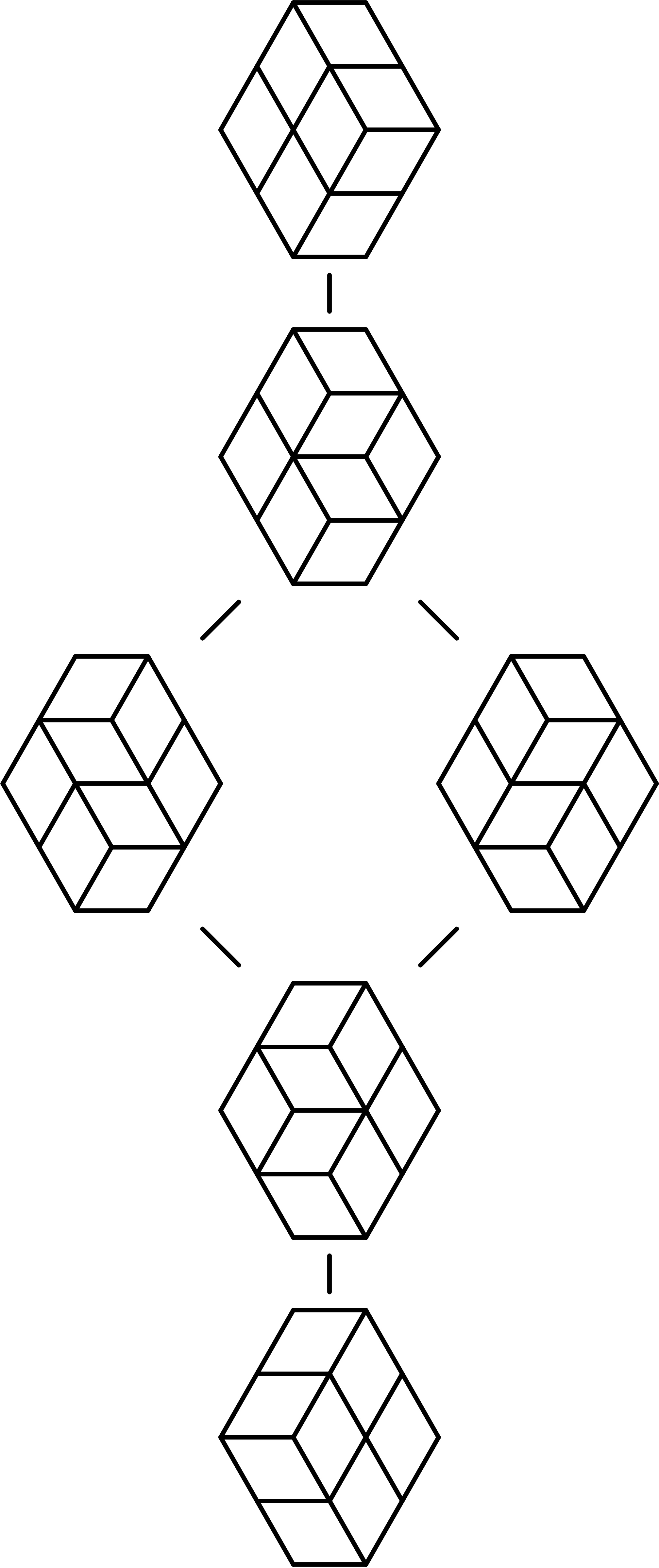}
\end{center}
\caption{Partial ordering of lozenge tilings.}
\end{figure}

This is a thinly disguised version of the
three-dimensional Young lattice ${\rm L}(a,b,c)$.
To see the connection,
imagine a collection of unit cubes nestling in the first octant of space,
so that in each of the three ``nestling directions''
$x \rightarrow -\infty$, $y \rightarrow -\infty$, and $z \rightarrow -\infty$,
a cube must be supported by either a wall
(one of the three coordinate planes)
or by another cube.
Viewing this collection from a point $(N,N,N)$ with $N$ suitably large,
the visible faces of the cubes,
together with the squares on the walls,
form a tiling of the plane by lozenges
(ignoring the scale-factor of $\sqrt{2/3}$).
If we restrict to collections of cubes
that fit inside an $a$-by-$b$-by-$c$ box,
then we may restrict the tiling of the plane
to a tiling of the hexagon with sides $a,b,c,a,b,c$
without losing any information.
The operation of twisting down an alternating cycle
corresponds to pushing down a piece of visible surface,
or equivalently,
to removing a cube from the collection.
These collections of cubes correspond to plane partitions
in an $a$-by-$b$ rectangle in which all parts are of size $c$ of less,
or equivalently, order-ideals in the poset defined
as the product of chains of cardinalities $a$, $b$, and $c$.
\end{example}

The ordering discussed in Example \ref{E21}
was first described by Thurston~\cite{Thur}.
The trick of turning three-dimensional Young diagrams
into matchings of a graph
was put to excellent use by Kuperberg~\cite{Kupe},
who used it in order to enumerate a hitherto intractable symmetry class
of plane partitions.

\begin{example} \label{E22}
Consider a plane region composed of unit squares
that can be tiled by {\bf dominoes}
(1-by-2 and 2-by-1 rectangles).
A basic move that turns one domino tiling into another
is the operation
that takes a single 2-by-2 block formed by two dominoes
and rotates it by 90 degrees.
Theorem \ref{thm2} implies that if the plane region under consideration
is simply connected, then it is possible
to get from any domino tiling of the region to any other
by means of such moves.
In the case where the region being tiled is a special sort of shape
called an Aztec diamond, the partial ordering on the set of tilings 
has an especially nice structure;~\cite{EKLP} gives details.
\end{example}

A concrete way of describing the height-function for domino tilings
as a function from the vertex-set to the reals
is to impose a checkerboard coloring on the squares
and to imagine an ant walking from vertex to vertex along edges
belonging to tile-boundaries;
the height increases by $1/4$ when the ant has a black square on its left
and decreases by $1/4$ when it has a white square on its left.
(This convention is related to that of~\cite{EKLP}
by a scale-factor of $1/4$ and is related to the convention 
of~\cite{Thur} by a scale-factor of $-1/4$.)
This picture for domino tilings, like the analogous picture 
for lozenge tilings, seems to have first been developed by Thurston.
Thurston's paper and the paper by Conway and Lagarias 
that inspired it~\cite{CoLa} deal only with the case in which
the graph $G$ underlying a set of tilings is very regular --
more precisely, the case in which
it is the intersection of the (planar) Cayley graph of a group
with a simply-connected region in the plane.
A good deal of the motivation for the theorems in this paper
was the belief that a more general formulation was possible and appropriate.

In the case of 1-factors of a graph $G$,
a very handy way of working with the lattice $\cM$
is to superimpose each matching $M$
with the fixed matching $M_{\hat{0}}$,
orienting the edges of $M$ from black to white
and the edges of $M_{\hat{0}}$ from white to black.
In this way, one obtains a spanning set of cycles,
where cycles of length 2 correspond
to edges that are common to $M$ and $M_{\hat{0}}$.
Let us omit the cycles of length 2.
One can then read this picture as a relief-map
in which the altitude changes by $+1$ or $-1$
as one crosses over a contour line
according to whether the contour-line
is directed toward the right or toward the left.
See Figure 11, in which $v^*$ is the center vertex on the top row
and $M_{\hat{0}}$ is the matching on the bottom of the figure.

\begin{figure}
\begin{center}
\includegraphics[width=0.5\textwidth]{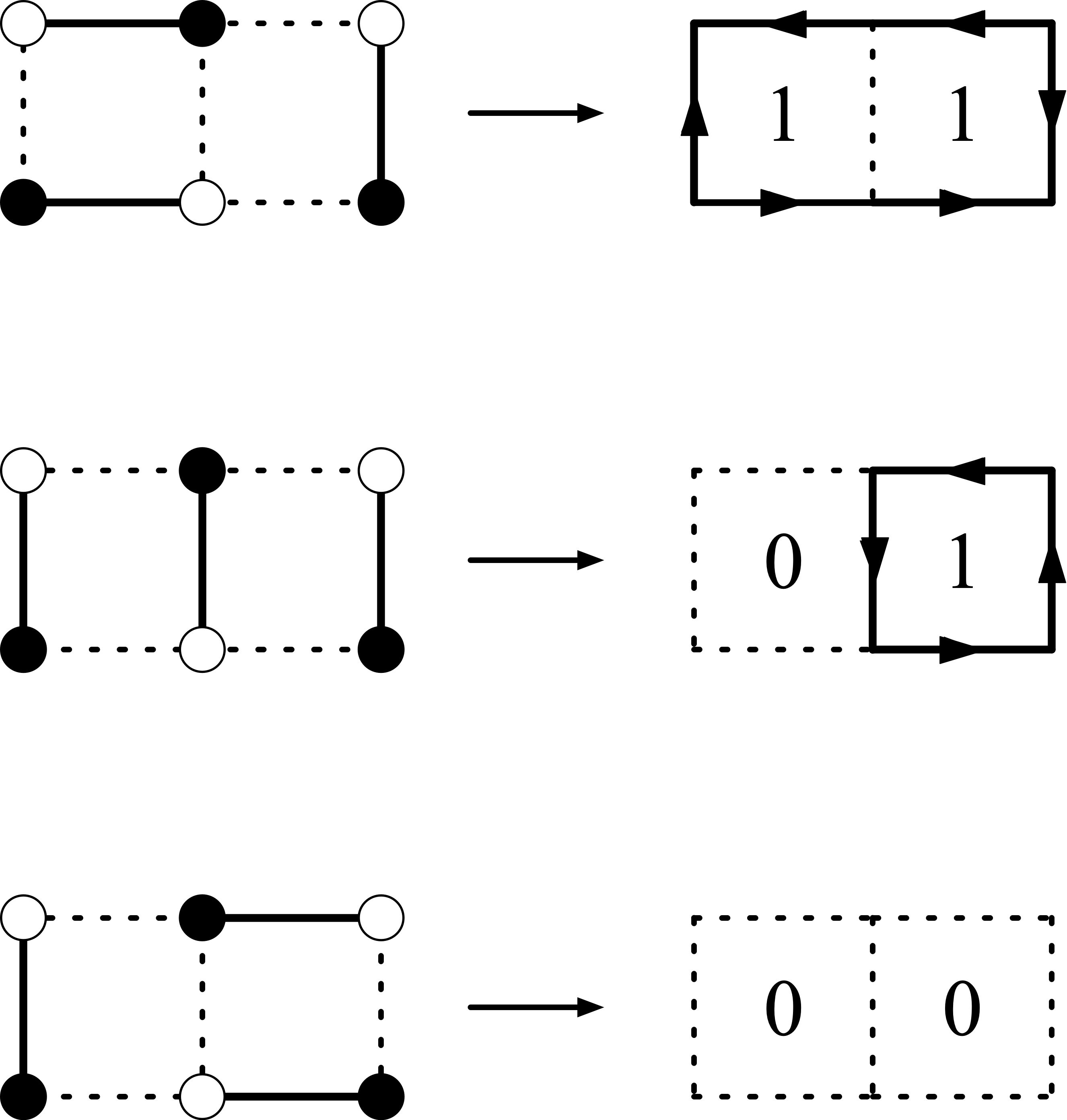}
\end{center}
\caption{Contour lines for matchings.}
\end{figure}

So far, we have been operating under the assumption that
every edge of $G$ belongs to some $d$-factors of $G$
but not to others.
We now indicate why this assumption is necessary.

\begin{example} \label{E23}
Consider the graph $G$ shown in Figure 12,
devised by William Jockusch.
This graph has twenty edges,
but only twelve of them
can participate in a matching of $G$;
the other eight have been shown as dashed edges.
The dark edges form three squares,
and every matching of $G$
is obtained by taking a matching of each square separately.
Thus, there are four matchings of $G$
in which the square of intermediate size is matched one way
and four matchings in which it is matched the other way.
It is easily seen that
the first four matchings are inaccessible from the other four
by means of face twists.
If we represent matchings of $G$
by $c$-orientations of the dual graph $G^\perp$,
we see that the dashed edges correspond
to two forced 4-cycles in $G^\perp$,
separating the square of intermediate size
from both the larger and the smaller squares.
Theorem \ref{thm1} is applicable,
but we must push down an entire cycle,
not just one vertex of a cycle.
Such a move corresponds to performing twist moves
on two squares simultaneously.
Equivalently, we could remove the dashed lines,
obtaining a graph on the sphere
with two simply-connected (4-sided) faces
and two non-simply-connected (8-sided) faces;
if we extend the notion of face-twists
to handle non-simply-connected faces
in the natural way,
we obtain the same partial ordering on the set of matchings.
\end{example}

\begin{figure}
\begin{center}
\includegraphics[width=0.4\textwidth]{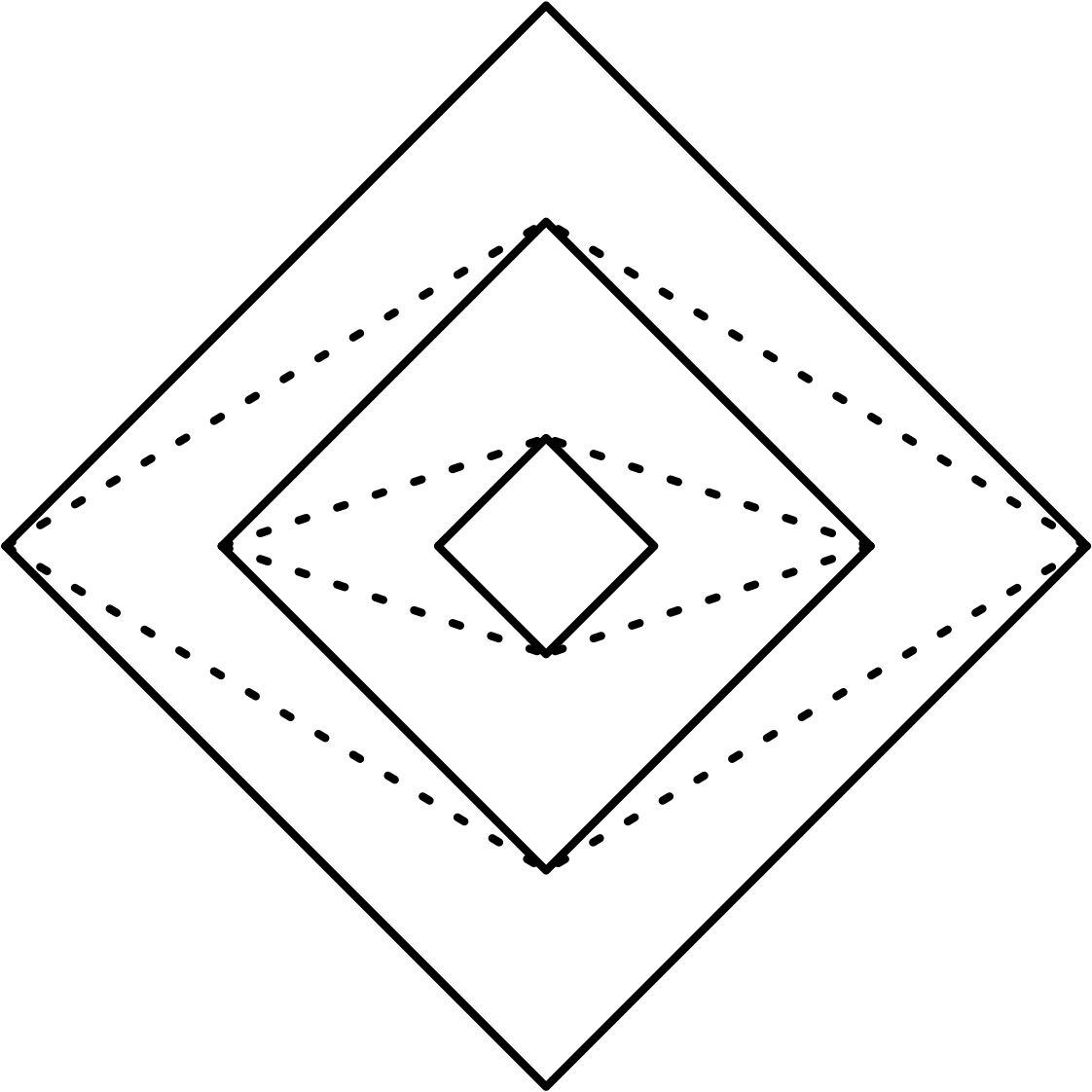}
\end{center}
\caption{Jockusch's example.}
\end{figure}

Now suppose $G$ is a connected bipartite graph on the torus.
Assume that every face of $G$
is a subset of a contractible patch of the torus.
Since the surface is orientable,
we may define twisting up and twisting down as before.
It is no longer the case
that all $d$-factors of $G$ are obtainable from one another
by twisting moves.
To see why,
recall that in the proof of Theorem \ref{thm2},
we showed that the orientation $R$
derived from a $d$-factor $M$
has fixed circulation around each elementary cycle $C$.
From this we concluded
that the circulation around {\it every} cycle $C$ was determined.
However, that conclusion follows only if
every cycle bounds.
This holds in the plane or on the sphere,
but not on the torus.
To specify the circulation around every cycle on the torus,
we must specify two additional parameters,
corresponding to the circulation around
non-bounding cycles that generate the homology of the torus.
We say that these two additional parameters $s$ and $t$
constitute the {\bf cohomology} of the $d$-factor,
and we represent it in the plane
by the point $(s,t)$.
(This representation has some geometrical dependence
on one's choice of two grid-paths that generate the homology of the torus,
but the properties of cohomology that we will illuminate via the 
two-dimensional picture are easily shown to be independent of these choices.)
See~\cite{STCR} for related earlier work on the cohomology of domino tilings.

\begin{example} \label{E24}
Let $G^\perp$ be a 4-by-4 grid on a torus
(shown as a grid on a 4-by-4 square
whose opposite edges are to be identified in pairs).
If we examine the domino tiling shown in
the left panel of Figure 13,
and try to use the ``ant-walk'' scheme
discussed in the paragraph following Example \ref{E22},
we get a globally-consistent height-function.
(The values along two of the four edges are shown.)
However, if we try to apply the ant-walk scheme
to the tiling shown in the right panel of Figure 13,
we find that the results we get are globally inconsistent;
as we travel all the way around the torus from bottom to top,
the ``height'' changes by 1.
These two tilings have cohomology $(0,0)$ and $(0,1)$, respectively.
Figure 14 shows the thirteen different cohomologies
that can arise from a domino tiling of the 4-by-4 torus.
\end{example}

\begin{figure}
\begin{center}
\includegraphics[width=0.7\textwidth]{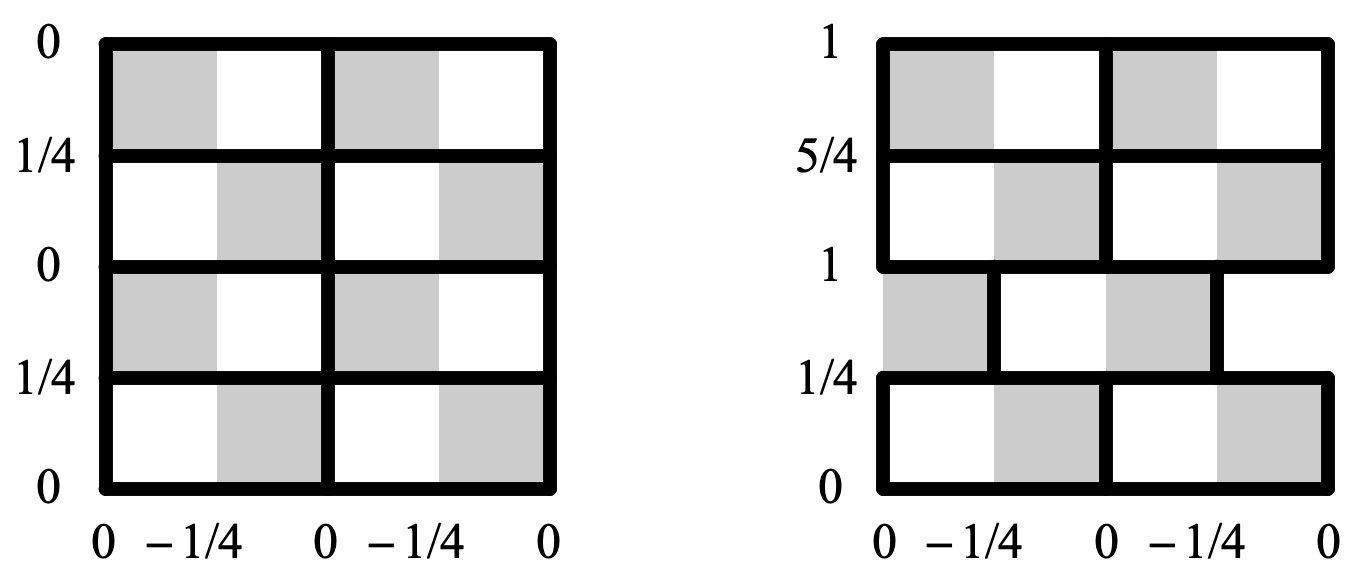}
\end{center}
\caption{Two non-cohomologous domino tilings.}
\end{figure}

\begin{figure}
\begin{center}
\includegraphics[width=0.4\textwidth]{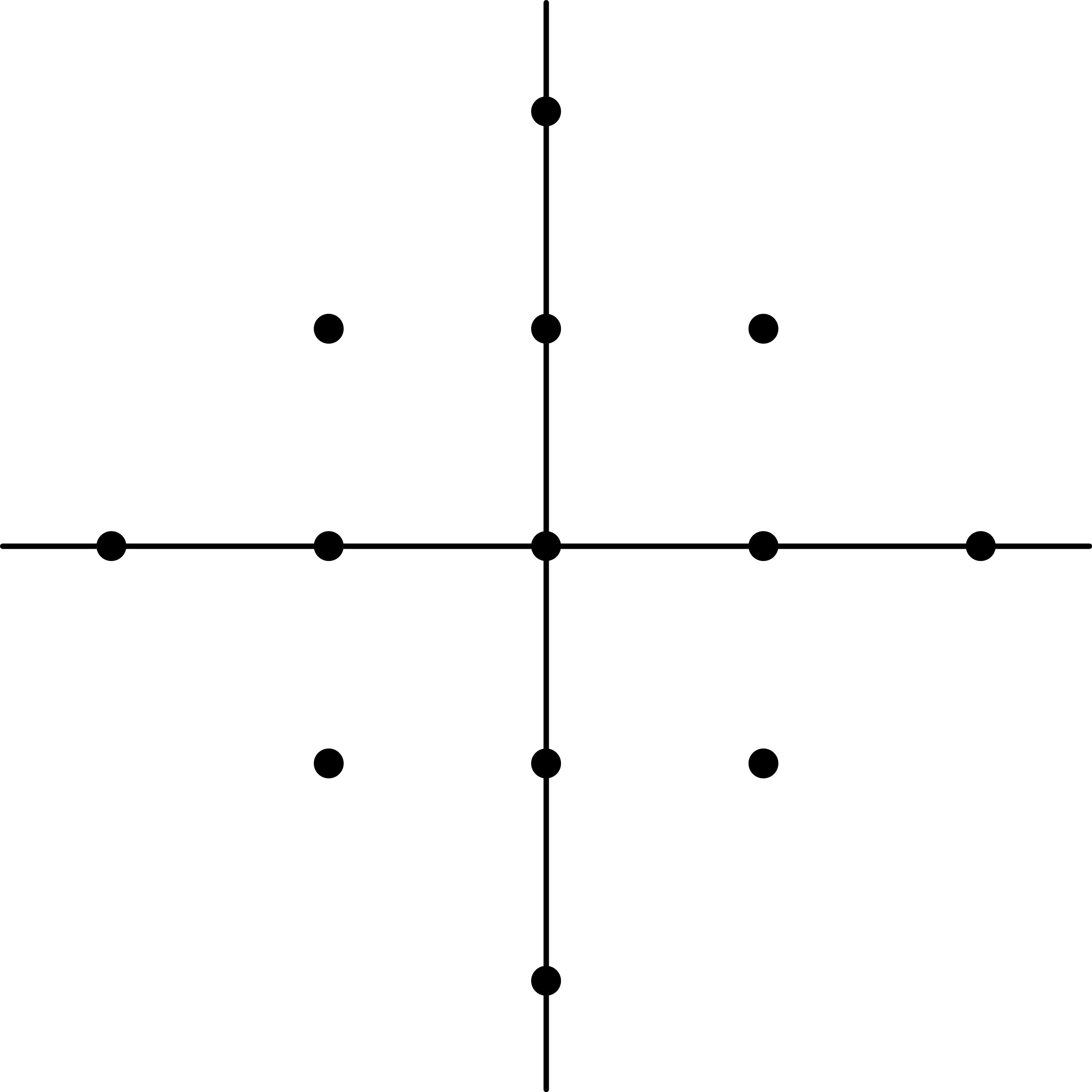}
\end{center}
\caption{Phase diagram for domino tilings of the 4-by-4 torus.}
\end{figure}

\begin{proposition} \label{PO}
If two $d$-factors of a torus graph $G$
have different cohomologies,
then neither one can be obtained from the other
by {\it any} form of local operation
that affects only the edges lying in
some contractible patch on the torus.
\end{proposition}

\begin{proof}
Consider two $d$-factors $M$ and $N$
obtained from one another by such a local operation.
Since the patch is contractible,
we can find paths
that generate the homology of the torus
but that avoid the patch.
Since $M$ and $N$ agree outside of the patch,
and since the behavior of the edges outside of the patch
suffices to determine cohomology,
we see that $M$ and $N$ must have the same cohomology.
\end{proof}

The inverse of Proposition \ref{PO} is almost true;
that is, when two $d$-factors have the same cohomology,
it is usually possible to get from one to the other,
not just using local operations of some unspecified type
but more specifically using face twists.
However, there are some exceptions.
For example, Figure 15 shows two domino tilings of the 4-by-4 torus
that have the same cohomology
but are not obtainable from one another
via face twists or any other sort of local move.

\begin{figure}
\begin{center}
\vspace{0.3in}
\includegraphics[width=0.5\textwidth]{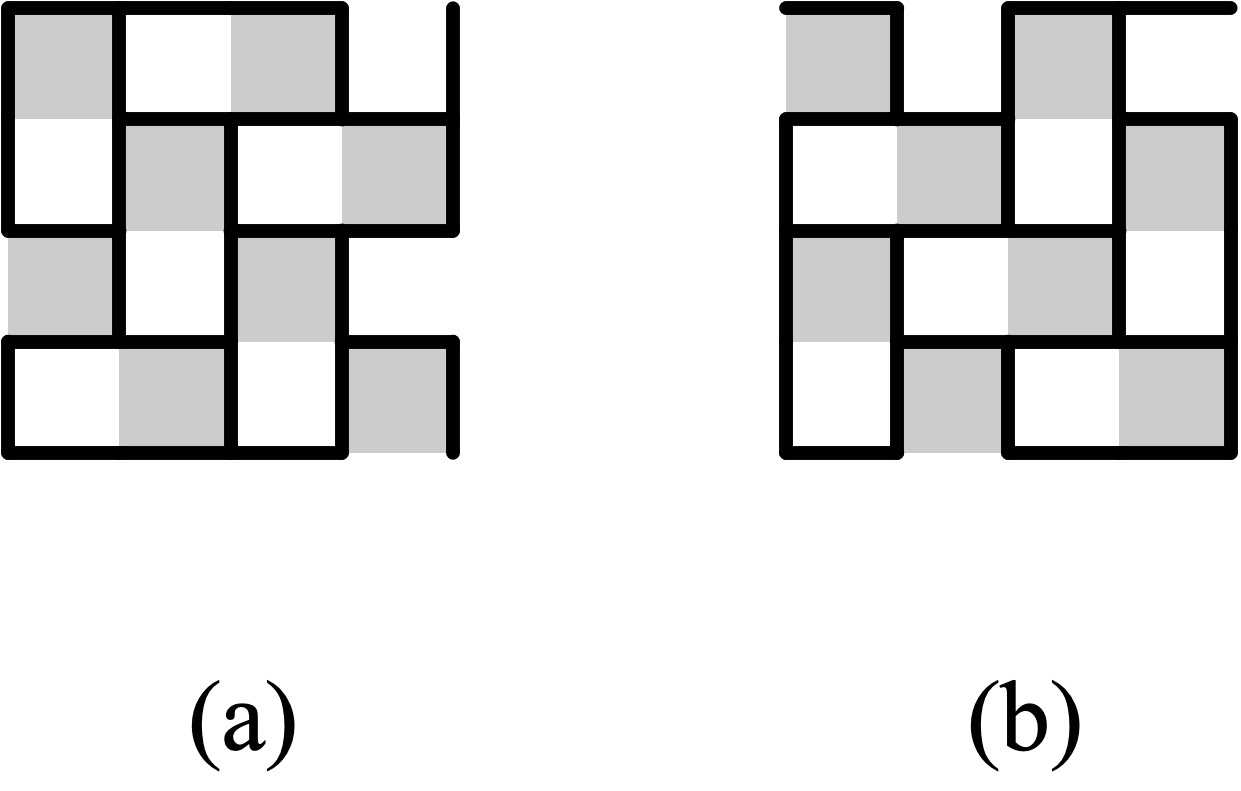}
\end{center}
\caption{Extremal cohomology.}
\end{figure}

Define the {\bf phase diagram} of a bipartite torus graph $G$
(relative to some degree-specification $d$)
as the set of ordered pairs arising
as cohomologies of $d$-factors of $G$,
or as the representation of this set in the plane.
(In the example just considered,
where $G$ and its dual $G^\perp$ are both isomorphic to
the 4-by-4 torus graph,
the phase diagram is given by Figure 14.)
We call a cohomology class of $d$-factors {\bf extremal}
if it is on the boundary of the convex hull
of the achievable cohomologies in the phase diagram.
The tilings shown in Figure 15
have the same, extremal cohomology.

\begin{proposition} \label{PP}
Let $M$ be a $d$-factor of the graph $G$.
If there exists a non-contractible forward cycle in
the associated orientation $R_M$ of $G^\perp$,
then $M$ has extremal cohomology.
\end{proposition}

\begin{proof}
If $R_M$ has a non-contractible forward cycle $C$ of length $m$,
then $R_M$ has circulation $m$ around $C$,
and no orientation of $G^\perp$ can have circulation $>m$ around $C$
as its length is only $m$.
Since the circulation around $C$
is an affine function of the cohomology,
we see that the cohomology of $R_M$
lies on a support line
for the set of achievable cohomologies,
and hence is extremal.
\end{proof}

\begin{proposition} \label{PQ}
Any $d$-factor $M$ of the graph $G$ having non-extremal cohomology
can be obtained via face twists
from any other $d$-factor having the same cohomology.
\end{proposition}

\begin{proof}
Assume $M$ has non-extremal cohomology.
Theorem \ref{thm1} implies
that any two $c$-orientations of $G^\perp$
can be obtained from one another by pushing down accessibility classes,
where $c$ is the circulation of $R_M$.
But Proposition \ref{PP} tells us that $c$ is acyclic,
so that accessibility classes are just vertices.
Since pushing down a vertex corresponds to
twisting down a face,
we are done.
\end{proof}

This theory can be extended to graphs
on orientable surfaces of higher genus $g$,
by way of $2g$-dimensional phase diagrams,
whose dimensions correspond to the $2g$ homology classes.

\section{Spanning Trees}

For the reader's convenience, I restate Theorem \ref{thm3}.

\bigskip

\noindent
{\bf Theorem \ref{thm3}:} {\em Let $\cT$ be the set of spanning trees
of a finite connected graph $G$ drawn on the sphere,
with $f^*$ a face of $G$
and $v^*$ a vertex of $G$ incident with $f^*$.
If we say that one tree $T \in \cT$ covers another tree $T' \in \cT$
exactly when $T'$ is obtained from $T$ by swinging down,
then the covering relation makes $\cT$ into
a distributive lattice.}

Before proving Theorem \ref{thm3},
let us reformulate it so that
the elements of our poset are not
spanning trees of $G$
but certain pairs consisting of a spanning arborescence of $G$
and a spanning arborescence of its dual graph.
Given a spanning tree $T$ of $G$,
there is a unique way to orient each edge
so that it points through the tree
towards the vertex $v^*$.
This oriented graph $\vec{A}$
is a {\bf spanning arborescence} rooted at $v^*$;
that is,
for every vertex $v$ there is a unique forward path
from $v$ to $v^*$ in $\vec{A}$.
In the other direction, we see that
if $\vec{A}$ is a spanning arborescence rooted at $v^*$,
the undirected edges associated with the directed edges of $\vec{A}$
constitute a spanning tree $T$.
Thus the spanning trees of $G$
are in 1-1 correspondence with
the spanning arborescences of $G$ rooted at $v^*$.

Let $G^\perp$ be the dual of $G$,
so that each edge $e$ of $G$
is associated with an edge $e^\perp$ of $G^\perp$.
Given a spanning tree $T$ of $G$,
the set of dual edges $\{e^\perp: e \not\in T\}$
is a spanning tree of $G^\perp$, called
the {\bf dual} of $T$.
We let $T^\perp$ denote the dual tree,
and $\vec{A}^\perp$ denote the unique orientation of $T^\perp$
that makes it a spanning arborescence rooted at ${f^*}^\perp$
(the vertex of $G^\perp$ dual to the face $f^*$ of $G$).
Thus the map
$T \mapsto (\vec{A}, \vec{A}^\perp)$
gives a bijection between spanning trees of $G$
and ordered pairs consisting of
a spanning arborescence $\vec{A}$ of $G$ rooted at $v^*$
and a spanning arborescence $\vec{A}^\perp$ of $G^\perp$ rooted at ${f^*}^\perp$
such that the spanning trees underlying $\vec{A}$ and $\vec{A}^\perp$
are dual to one another.

Using the pair $(\vec{A},\vec{A}^\perp)$,
we may reformulate our original definition of swinging down.
Suppose that $v$ is a vertex of $G$
and $e,e'$ are edges of $G$
such that $e'$ is the clockwise successor of $e$ at $v$.
Let $\vec{e}$ and $\vec{e}\,'$ be
the respective orientations of $e$ and $e'$ away from $v$,
and let $\vec{e}^\perp$ and ${\vec{e}\,'}{^\perp}$
be the respective dual edges oriented so as to point away from $f^\perp$,
where $f$ is the face of $G$
that contains angle $eve'$.
(See Figure 16.)
Let $T$ be a spanning tree of $G$,
with $(\vec{A},\vec{A}^\perp)$
the associated pair of arborescences.

\begin{figure}
\begin{center}
\includegraphics[width=0.4\textwidth]{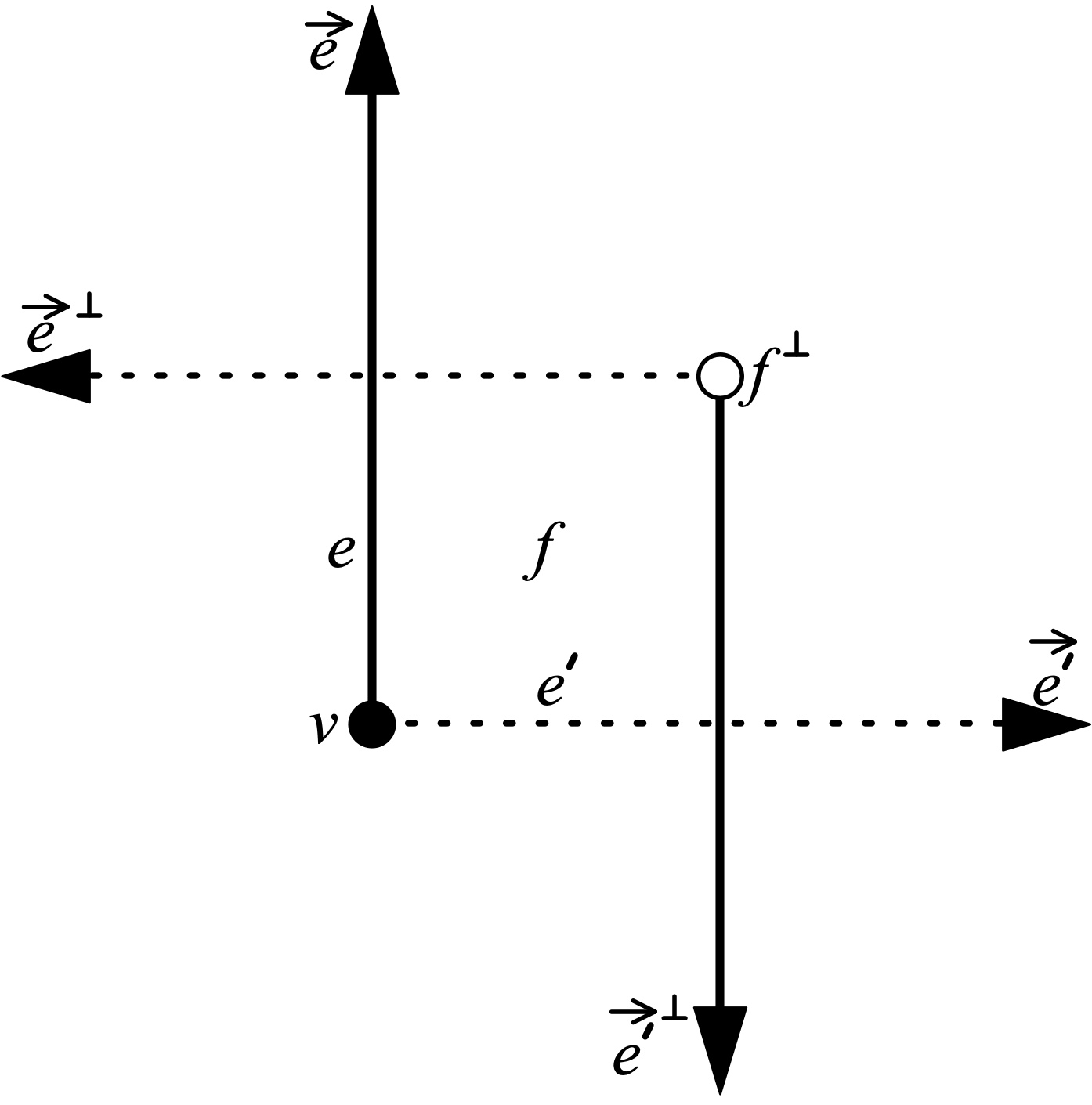}
\end{center}
\caption{Pivoting dual arborescences.}
\end{figure}

\begin{proposition} \label{PR}
Under the foregoing assumptions,
the angle $eve'$ is positively pivotal
if and only if the following conditions are satisfied:
\begin{enumerate}
\item[$(1')$] $\vec{e} \in \vec{A}$ and $\vec{e}\,' \not\in \vec{A}$;
\item[$(2')$] the symmetric difference
$\vec{A}' = \vec{A} \symmdiff \{\vec{e},\vec{e}\,'\}$ is
a spanning arborescence of $G$ rooted at $v^*$; and
\item[$(3')$] the symmetric difference $(\vec{A}^\perp)' =
\vec{A}^\perp \symmdiff \{\vec{e}^\perp,{\vec{e}\,'}{^\perp}\}$ is
the spanning arborescence of $G^\perp$ rooted at ${f^*}^\perp$
that is dual to $\vec{A}'$.
\end{enumerate}
\end{proposition}

\begin{proof}
Referring back to the definition of ``positively pivotal''
given in Section 1, we see that
(1) is encoded by $(1')$;
(2) is encoded by $(2')$;
(3) and (4) are encoded by our definition of $\vec{A}$ and $\vec{A}'$
as spanning arborescences with specified roots; and
(5) is encoded by $(3')$.
\end{proof}

\begin{proposition} \label{PRa}
Under the same assumptions as prevailed in Proposition \ref{PR},
the angle $eve'$ is positively pivotal
if and only if
$\vec{e} \in \vec{A}$,
$\vec{e}\,' \not\in \vec{A}$,
$\vec{e}^\perp \in \vec{A}^\perp$,
${\vec{e}\,'}{^\perp} \not\in \vec{A}^\perp$.
\end{proposition}

\begin{proof}
Easy.
\end{proof}

Using this picture of swinging down we can now prove
that swinging down at a single vertex cannot be repeated indefinitely.

\begin{proposition} \label{PS}
If $T$ is a spanning tree of $G$, it is impossible to swing down an edge 
$n$ times in succession around a vertex $v$, where $n = \deg v$.
\end{proposition}

\begin{proof}
Let $C^\perp$ be the simple cycle in $G^\perp$ that encircles the vertex $v$.
The spanning arborescence $\vec{A}^\perp$ of $G^\perp$
associated with $T$ contains $n-1$ directed edges of $C^\perp$.
Each time we swing down, we replace a clockwise edge of $C^\perp$
by a counterclockwise edge, so we can do it at most $n-1$ times.
\end{proof}

(Of course, it may become possible to swing down
many more times around the vertex $v$
if we intersperse these swinging-down moves
with swinging-down moves at other vertices.)

We can now prove Theorem \ref{thm3}.

\begin{proof}
The main idea is a generalization of
a trick due to Temperley~\cite{Temp} and discussed by Lov{\'a}sz~\cite{Lova}
(problem 4.30, pages 34, 104, 243--244);
subsequently, Burton and Pemantle~\cite{BuPe} generalized Temperley's idea
and applied it to infinite graphs drawn in the plane.
Let $\cH(G)$ be the Hasse diagram of the face-poset of $G$,
viewed as a graph.  Concretely, $\cH(G)$ is a graph
with a node $\overline{v}$ corresponding to each vertex $v$ of $G$,
a node $\overline{e}$ corresponding to each edge $e$ of $G$,
and a node $\overline{f}$ corresponding to each face $f$ of $G$,
with an edge joining two nodes in $\cH(G)$
if the corresponding elements in $G$ are either
an edge and one of its endpoints
or an edge and one of the faces it bounds.
For instance, if $G$ is the graph shown in Figure 17(a),
$\cH(G)$ is the graph shown in Figure 17(b).
We can view $\cH(G)$ as the result
of jointly embedding $G$ and $G^\perp$ on the sphere.
Let $H=H(v^*,f^*)$ be the graph obtained from $\cH(G)$
by deleting the nodes $\overline{v^*}$ and $\overline{f^*}$
(along with all incident edges in $\cH(G)$).
Euler's formula $|V|-|E|+|F|=2$ implies that $|V|+|F|-2=|E|$, so that
$H$ is a balanced bipartite graph.  I will give a bijection between
the spanning trees of $G$ and the matchings of the bipartite plane graph $H$.

\begin{figure}
\begin{center}
\includegraphics[width=4in]{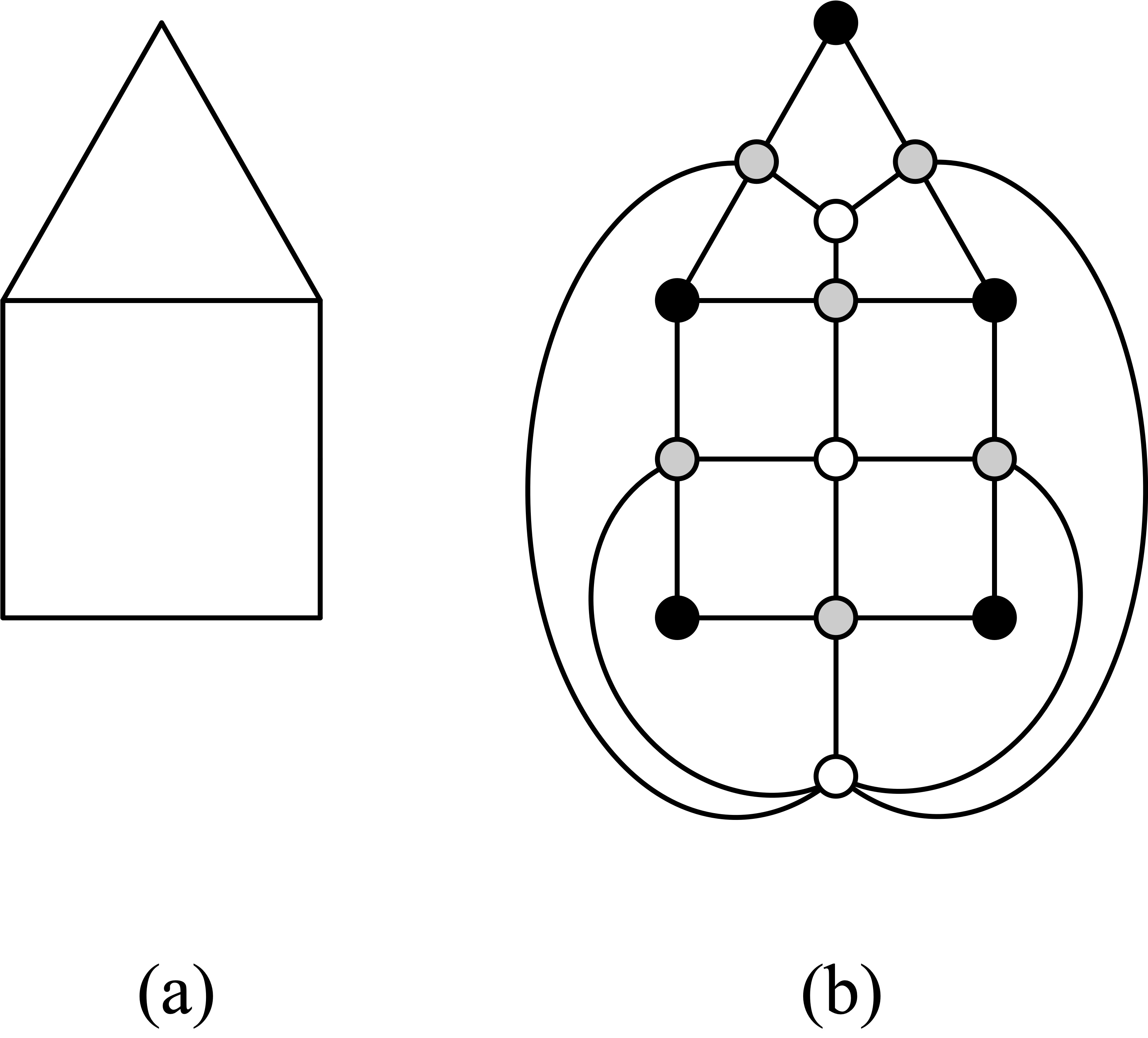}
\end{center}
\caption{The Hasse diagram graph.}
\end{figure}

Every spanning tree $T$ of $G$
(together with the associated pair $(\vec{A}, \vec{A}^\perp)$)
determines a matching $M$ of $H$ as in Figure 18.
Specifically, for each vertex $v$ of $G$,
pair $\overline{v}$ with the unique $\overline{e}$
such that $v$ is incident with $e$
and $e$ is pointing away from $v$ in $\vec{A}$,
and  pair $\overline{f}$ with the unique $\overline{e}$
such that $e$ is incident with $f$
and $e^\perp$ is pointing away from $f$ in $\vec{A}^\perp$.

\begin{figure}
\begin{center}
\includegraphics[width=0.6\textwidth]{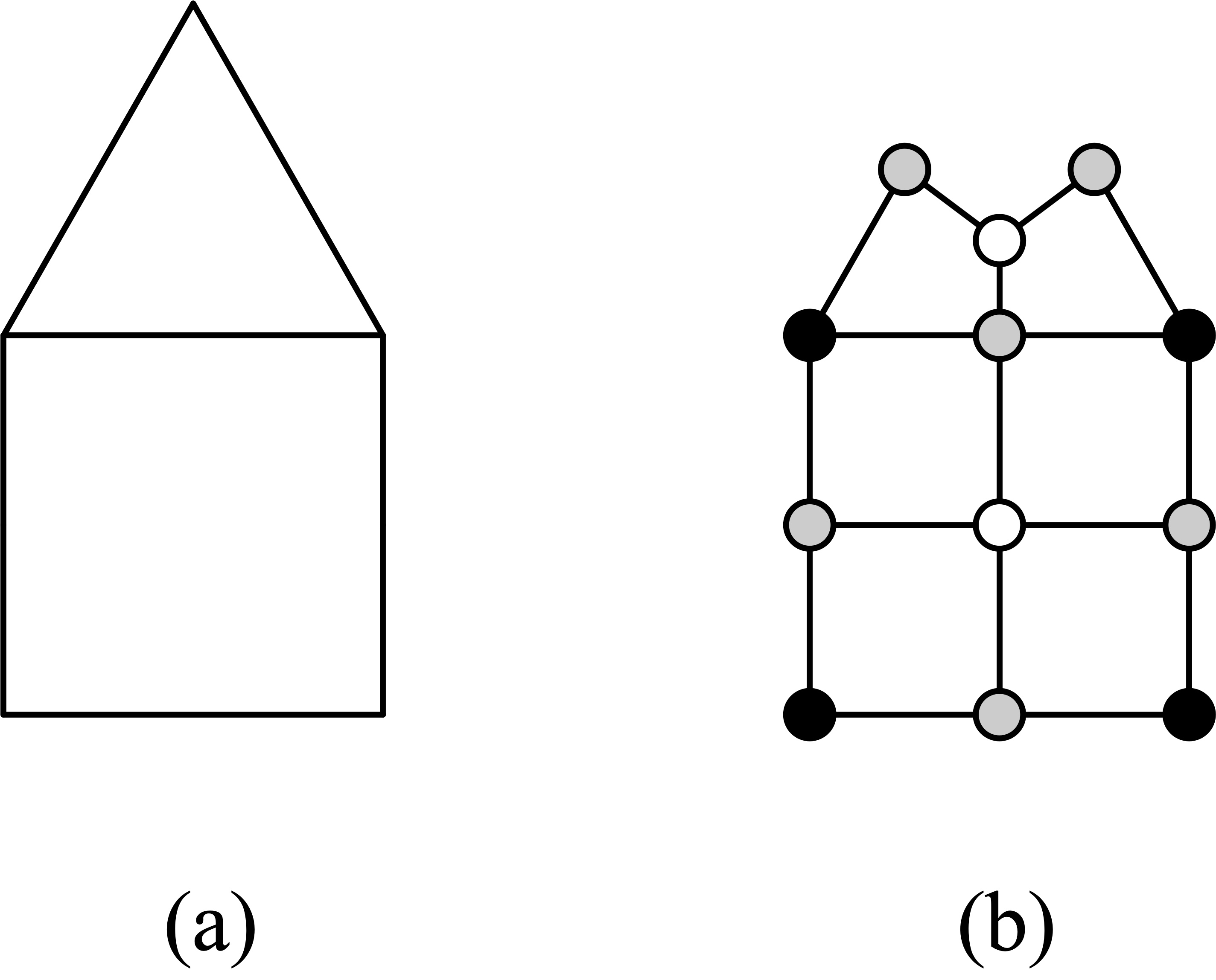}
\end{center}
\caption{From spanning trees to matchings.}
\end{figure}

To verify that this gives a matching $M$ of $H$,
it suffices to show that no edge-node $\overline{e}$ is paired twice.
But this could only happen if we had $e \in T$ and $e^\perp \in T^\perp$,
contradicting the definition of $T^\perp$.
From the matching $M$ we can easily recover $T$ as $\tilde{T}_M$,
defined as the set of edges $e$ such that
$\overline{e}$ is paired with a vertex-node in $H$
under the matching $M$.
Hence the mapping that turns $T$ into $M$ is injective.

To show that the mapping is surjective as well,
let $M$ be a matching of $H$, with $\tilde{T}_M$ as above.
We must show that $\tilde{T}$ is a spanning tree.
Since $\tilde{T}$ spans $G$ and has $|V|-1$ edges,
it suffices to prove that $\tilde{T}$ is acyclic.

Suppose $\tilde{T}$ contains a cycle $C$, say of length $n$.
$C$ divides the sphere into two (open) parts, or hemispheres,
one of which contains both $v^*$ and $f^*$
and the other of which contains neither.
I claim that each hemisphere contains an odd number of nodes of $\cH(G)$
and hence an odd number of nodes of $H$ as well.
For, suppose we modify $G$ by replacing either of the two hemispheres
by a single face.
By Euler's formula, the quantity $|V|+|E|+|F|$
(the number of ``elements'' of $G$) in the resulting graph must be even.
Since there are an even number of elements of $G$ on the cycle $C$
($n$ vertices and $n$ edges) and an odd number in the modified hemisphere 
(1 face), the unmodified hemisphere must have an odd number of elements 
of $G$ as well.

Since the edges of $C$ disconnect $H$ into parts lying in the two hemispheres,
$M$ must match each hemisphere within itself.
But this is impossible, since each hemisphere has been shown to contain
an odd number of nodes of $H$.  Hence no such cycle $C$ can exist.

This completes the validation of the bijection between spanning trees of $G$
and matchings of $H$.  Since $H$ is a connected bipartite graph
drawn on the plane, we can apply Theorem \ref{thm2}.
It is easily checked that modifying a spanning tree $T$ of $G$
by swinging down a directed edge at $v$ through angle $eve'$
corresponds to modifying the associated matching $M$
by performing a twist on the face of $H$ with nodes
$\overline{v}$, $\overline{e}$, $\overline{f}$, and $\overline{e'}$,
where $f$ is the face belonging to angle $eve'$.
\end{proof}

Note that Theorem \ref{thm3} gives us as many as $2|E|$
distinct orderings of the set of spanning trees of $G$
(since that is the number of ways in which
we can choose a vertex $v^*$ and an adjacent face $f^*$).

It is possible to bypass the invocation of Theorem \ref{thm2}
and give a direct description of a height-function
for spanning trees of a graph $G$.
Such a height-function is a mapping from the angles of $G$ to the real numbers,
where an angle of $G$ is specified by a vertex and an adjacent face.
Every time we swing down an edge of spanning tree through a particular angle,
we decrease the height associated with that angle by 1.

One especially nice class of graphs $G$ to which Theorem \ref{thm3} applies
(suggested by Curtis Greene~\cite{Gree}) is a class of graphs
I call {\em Hamiltonian outerplanar graphs}.
These are graphs that are outerplanar
in the sense that they are drawn in the plane,
in such a way that all the vertices are adjacent 
to the unbounded face $f^*$; 
additionally they are Hamiltonian in that there is
a Hamiltonian path bordering the unbounded face.
In a speculative vein, we might also consider what I call
{\em Hamiltonian weakly outerplanar graphs},
in which edges of $G$ are allowed to cross
inside the Hamiltonian circuit.
More specifically, two edges of $G$ that have no endpoints in common
are permitted to cross at a single point
that is distinct from all the vertices of $G$
and lies on no other edge of $G$,
but two edges of $G$ that share an endpoint cannot have extra crossings.
See Figure 19. Although Theorem \ref{thm3} says nothing about this situation,
it appears (see the example presented at the end of this section)
that the conclusion of the theorem might apply
under weaker hypotheses.

\begin{figure}
\begin{center}
\includegraphics[width=0.2\textwidth]{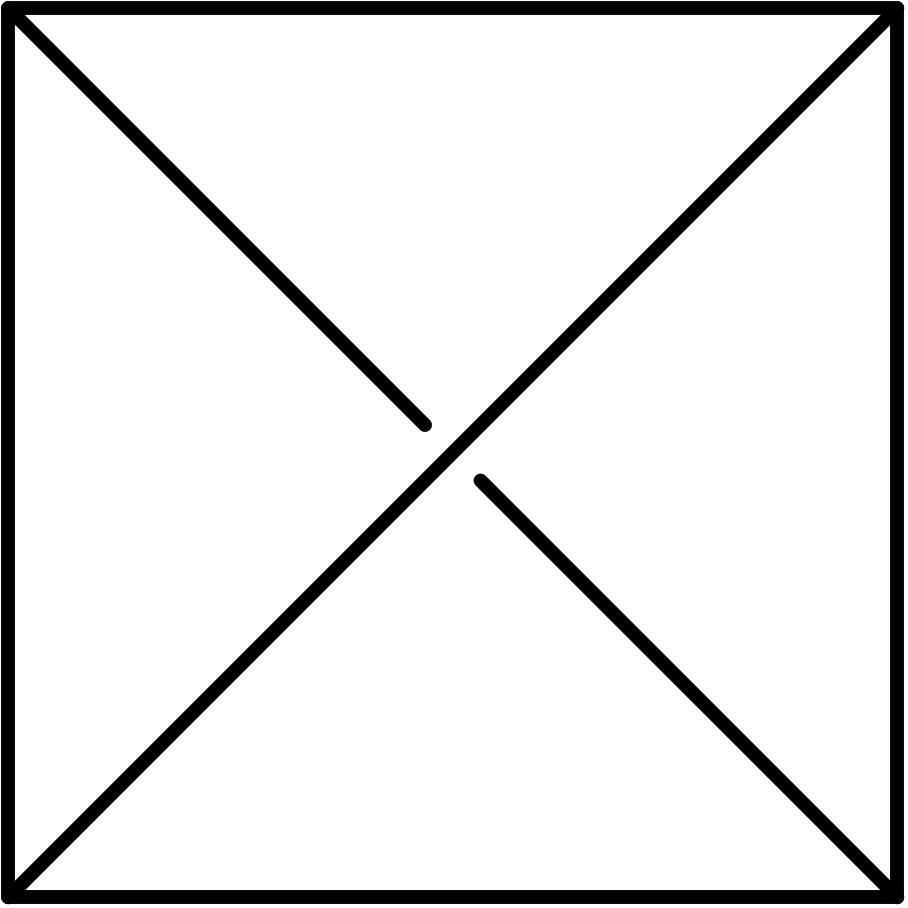}
\end{center}
\caption{A Hamiltonian weakly outerplanar graph.
The two diagonal edges cross but do not intersect.}
\end{figure}

The crossing-free case has many nice features.
The first is that condition (5) in the definition of a swinging-down move
follows from the first four:
for, the simple cycle in question will be a convex polygon
that contains the face $f$ but not the face $f^*$.
A second feature is that we can say exactly
what the minimum and maximum elements of the lattice are:
if we number the edges of the Hamiltonian cycle
as $e_1$, $e_2$, ..., $e_n$
in clockwise order,
where $e_1$ and $e_n$ have $v^*$ in common,
then the minimum element is the tree (a path, actually)
consisting of
all the edges in the cycle except $e_n$
and the maximum element is the tree consisting of
all the edges in the cycle except $e_1$.
A third feature is that for every vertex $v$,
the face $f^*$ blocks one from swinging down an edge all the way around $v$,
even if (unlike the situation of Proposition \ref{PS})
one intersperses swinging-down moves at $v$ with swinging-moves elsewhere.
That is, a given angle of a face of $G$
can be used for swinging-down only once.

These properties give us
a concrete interpretation of
the join-irreducibles of the lattice $\cT$
as angles of $G$,
as shown in Figure 20.
In the right panel of Figure 20,
the empty circles represent the angles of $G$
at vertices other than $v^*$,
and constitute a set $P$;
the arrows between the empty circles give a covering relation
that makes $P$ a poset
isomorphic to the poset of join-irreducibles of $\cT$.
The prescription for drawing the arrows is as follows:
First draw the dual of $G$, with the special vertex ${f^*}^\perp$
corresponding to the unbounded region $f^*$.
Erase all edges of this dual graph
that contain the vertex ${f^*}^\perp$
except for the edge dual to the edge of $G$
that joins $v^*$ with its counterclockwise neighbor.
This modified dual is a tree; make it into
an arborescence $\vec{A}$ rooted at ${f^*}^\perp$,
as in the left panel of Figure 20.
Now, at each vertex $v$ of $G$ other than $v^*$,
connect the circles (representing angles of $G$) by clockwise arrows,
except that there must be no arrow going through the unbounded region $f^*$.
Lastly, at the boundary of each face $f$ of $G$ other than $f^*$,
connect the circles by clockwise arrows, except that there must be
no arrow on the edge of $f$ that is dual to the directed edge of $\vec{A}$
that points away from $f^\perp$.

\begin{figure}
\begin{center}
\includegraphics[width=5in]{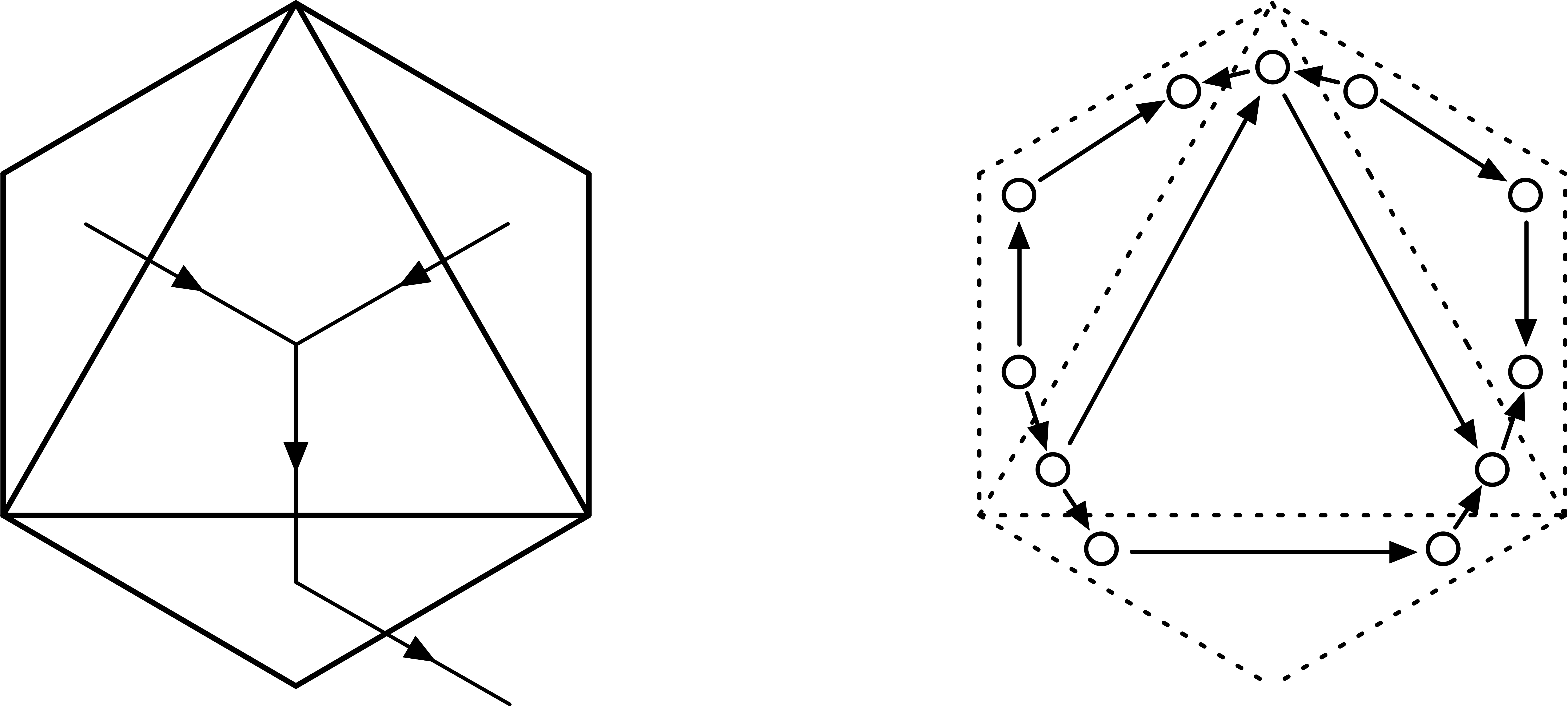}
\end{center}
\caption{The join-irreducibles for a Hamiltonian outerplanar graph.}
\end{figure}

Notice that the idea of swinging down an edge is perfectly well-defined
even when the outer-Hamiltonian graph $G$ has crossings,
since one still has a cyclic orientation for
the edges that emanate from a vertex.
In this context there does not appear to be a dual graph,
so the version of ``positively pivotal'' that I have in mind
retains conditions (1') and (2') from Proposition \ref{PR}
but drops condition (3').
For instance, let $G$ be the complete graph $K_n$,
drawn in the plane with crossings
(Figure 19 shows the case $n=4$),
let $v^*$ be some vertex of $G$,
and let $f^*$ be the external $n$-sided face.
Then one can define a partial ordering $\cT$
on the set of spanning trees of $G$,
along the lines of the statement of Theorem \ref{thm3}.
In 1990, Michelle Thompson proved~\cite{Thom} that
the rank generating function of $\cT$ is
$(1+q+q^2+...+q^{n-1})^{n-2}$,
which is a strengthening of Cayley's formula $n^{n-2}$
for the number of spanning trees of a complete graph~\cite{Lova}.

\section{Earlier and later related work}

The four preceding sections (with minor stylistic differences) appeared in 
a preprint released in 1993 that was posted on the arXiv in 1995
but not published until now.
In the past three decades a great deal of work has been done in this area,
much of which builds on that preprint. In this section I attempt a brief,
mostly chronological, and necessarily incomplete survey of that work.
First, however, I will acknowledge work that I was unaware of back in 1993
but which anticipated the height-function construction in various ways.

Note: I will make no attempt to distinguish between the date a result was 
proved, the date a paper proving the result was released, and the date 
the paper was published, so even when I associate a result with a
particular year, the year should be taken only as a rough approximation.
In no way should my attempt at an overview be construed as any assertion 
that one researcher or group of researchers proved a result before another, 
except where I explicitly mention influence.

In 1977, Henk van Beijeren~\cite{Beij} gave a height-function representation
for the square ice model. Under a natural bijection between square-ice states 
and alternating-sign matrices (see for instance~\cite{Prop}),
van Beijeren's representation is equivalent to the height-function for 
alternating-sign matrices presented earlier in the current paper.
In physicist terms, van Beijeren gave a ``solid-on-solid'' (SOS)
interpretation of the two-dimensional model; finding such interpretations
was an important tactic of researchers in statistical mechanics (see
for instance~\cite{ABF}). To see how pervasive the technique of introducing
height functions had become in statistical physics of lattice models, see 
for instance~\cite{MN} and~\cite{MNMS} and the references contained therein.

In 1982, Henke Bl\"ote and Henke Hilhorst~\cite{BH} came up with 
an SOS interpretation for dimers in a hexagonal lattice,
equivalent to the one presented in this article.
This interpretation was reinvented in 1984 by Veit Elser~\cite{Else}, 
who used the connection between hexagonal-lattice dimer-configurations 
and plane partitions to draw conclusions about how bulk entropy 
can be sensitive to the shape of the boundary of a region.
This boundary-dependence phenomenon was noticed independently 
a decade later by Horst Sachs and Holger Zernitz~\cite{SZ},
and an analogous observation for the square grid had been made
by Grensing, Carlsen, and Zapp~\cite{GCZ} as early as 1980;
however, those two papers did not discuss SOS interpretations
or anything resembling height-functions.

In 1983, knot-theorist Louis Kauffman~\cite{Kauf}, 
in constructing his state-sum model for the Alexander-Conway polynomial, 
devised a partial ordering related to knot diagrams
and showed in his Clock Theorem (Theorem 2.5) 
that the partial order gave rise to a lattice. 
Kauffman's analysis shows that this lattice is in fact distributive, 
though this fact was not stated explicitly 
until the 1986 article of Patrick Gilmer and Richard Litherland~\cite{GL}.
As Moshe Cohen and Mina Teicher~\cite{CT} would later explain, 
Kauffman's lattice is a special case of the distributive lattice 
on perfect matchings of a graph; Figure 3 of their article shows 
how Kauffman's clock move relates to the twist move for matchings.
The work of Gilmer and Litherland was later 
extended by Hine and K\'alm\'an~\cite{HK}.
See also recent work by Le and Mathews~\cite{LM}.

In 1986, Oliver Pretzel wrote the two articles~\cite{Pr1,Pr2} 
discussed earlier. It is worth stressing that his results 
do not require that the graph under consideration be planar.
His construction depends on a choice of a sink, accounting for the fact 
that my construction, in a dual setting, depends on a choice of a face.

In 1990, pioneering articles 
by John Conway and Jeffrey Lagarias~\cite{CoLa} 
and William Thurston~\cite{Thur} were published
(the latter having been influenced by a preprint of the former);
Thurston gave detailed descriptions of height-functions 
for both domino tilings and rhombus tilings,
and used height functions to devise a linear-time algorithm
to test a region for tileability and, if a tiling exists, produce one.
In the meantime, Leonid Levitov~\cite{Levi} came up with 
a height-function for the dimer model on a square grid 
that is equivalent to Thurston's height-function for domino tilings 
by way of the natural duality between dimers and dominoes.

The work of Conway, Lagarias, and Thurston~\cite{CoLa,Thur} 
inspired me to seek regions in the square grid 
whose boundary height-functions vis-a-vis domino tilings 
would resemble those of regular hexagons vis-a-vis lozenge tilings,
and thus led me to come up with Aztec diamonds. As it turned out, 
these regions had already been studied in the aforementioned work of
Grensing, Carlsen, and Zapp~\cite{GCZ}, who stated (in different words) 
that the number of domino tilings of the Aztec diamond of order $n$ is 
$2^{n(n+1)/2}$, but the authors of~\cite{EKLP}, who presented four proofs 
of the formula, did not learn of the work of Grensing et al.\ until later. 
In any case, the niceness of the formulas counting domino tilings 
of Aztec diamonds and lozenge tilings of hexagons convinced me that 
a more general understanding of height-functions for dimer models 
for plane graphs was worth developing, and motivated the 1993 preprint.

In 1992, Claire and Richard Kenyon~\cite{KK} used a notion of height 
based on the (usually nonabelian) groups of Conway and Lagarias; 
they defined height as the distance from an element of the group 
to the identity element in the Cayley graph. They used this kind of
height both to construct algorithms for determining tileability and 
to prove that in a number of tiling contexts, certain simple local 
moves, applied in sequence, suffice to turn any tiling into any other.

In 1993, working independently of the aforementioned
mathematicians, physicists, and chemists, the computer scientists
Samir Khuller, Joseph Naor, and Philip Klein~\cite{KNK} described 
a distributive lattice structure on the set of integer circulations 
in a planar graph. Their potential functions correspond to my height 
functions, but their work is more general in that they go beyond the
setting of $\alpha$-orientations.

In 1994, in yet another independent line of work, 
Patrice Ossona de Mendez~\cite{Osso}, showed that the set of 
indegree-constrained orientations has a distributive lattice structure. 
Subsequent work by Enno Brehm~\cite{Breh} extended this result
to 3-orientations.

In 1995, Nicolau Saldanha, Carlos Tomei, Mario Casarin Jr., and 
Domingos Romualdo~\cite{STCR}, drawing upon the work 
of Conway and Lagarias~\cite{CoLa} and Thurston~\cite{Thur}, 
came up with a generalization of height functions 
they called {\em height sections} that can be used when the region being tiled 
is not simply connected; in such cases the local change in height (defined in 
the most natural manner) can have non-zero sum along a non-contractible cycle, 
making it impossible to define a single-valued height function. 

One important application of putting a lattice structure 
on a large set of combinatorial objects is sampling.
In its full generality, the problem of sampling 
uniformly at random from a large set is intractable;
at the same time, approximately uniform random samples
can often be generated via Markov chain methods.
In 1995 Propp and Wilson~\cite{PW}
introduced the method of ``coupling from the past'' (CFTP),
one specific form of which (monotone CFTP)
could eliminate ``burn-in time'' as a source of discrepancy
between the target distribution and the sampling distribution.
The method was applied not just to ice and dimer models
but also to the random cluster model of Fortuin and Kasteleyn
and thence to the Ising model. For an exploration of the method
in the setting of $\alpha$-orientations of plane graphs, see~\cite{FH}.
Imposing (or discovering) a distributive lattice structure
on a large set isn't useful only for random sampling;
it can also expedite other sorts of operations on the set~\cite{HMNS}. 

In 1996, unaware of my (unpublished) work, Thomas Chaboud~\cite{Ch} 
used height functions (which he dubbed altitude functions) to extend 
Thurston's algorithm to determine whether a given regular bipartite 
planar graph has a perfect matching, and if so, construct one.

In 1998, \'Eric R\'emila~\cite{Re1} used height functions on Cayley 
graphs to extend Thurston's method to tilings in the triangular grid 
in which the tiles are ``leaning bars'' and equilateral triangles.

In 2001, Joakim Linde, Cristopher Moore, and Mats Nordahl~\cite{LMN} 
defined height functions for higher-dimensional analogues
of rhombus tilings and used coupling-from-the-past
to explore the behavior of random tilings.

In 2002, Matthieu Latapy and Cl\'emence Magnien~\cite{LaMa}
proved a kind of converse to my Theorem \ref{thm1}:
just as the constrained orientations of a finite graph can
always be modeled as elements of a finite distributive lattice,
every finite distributive lattice can be modeled as the set
of constrained orientations of a finite graph.

Also in 2002, Scott Sheffield~\cite{Shef} devised 
a multidimensional height function for ribbon tilings and used them to prove
that all $n$-ribbon tilings of a simply-connected polyomino
can be obtained from each other by a sequence of mutations of the simplest kind
(``2-flips'' that replace two ribbon tiles by two other ribbon tiles),
as had been conjectured by Igor Pak~\cite{Pak1}.

In that same year, Cristopher Moore, Ivan Rapaport 
and \'Eric R\'emila~\cite{MRR2} 
defined a height function and proved a local connectivity property for 
certain sets of Wang tiles (colored square tiles with color-based adjacency 
constraints).

In 2003, Oliver Pretzel, in response to my preprint, 
published his own proof~\cite{Pr3} of my Theorem \ref{thm1},
along with a discussion of how the lattice changes when the chosen sink 
varies. In that same year Peter Che Bor Lam and Heping Zhang
published a paper~\cite{LaZh} whose main content was Theorem \ref{thm2},
which they had discovered independently and proved in a different way.

Also in 2003, Igor Pak published a survey article~\cite{Pak2} in which 
he articulated, among other things, a proposed general framework for 
height functions.

In 2004, \'Eric R\'emila~\cite{Re2} provided his own proof 
(along lines similar to mine) in the special case of domino tilings 
and made a more detailed study of the join-irreducibles of the lattice,
relating them to the geometry of the region being tiled.

Also in 2004, Stefan Felsner~\cite{Felsner}, 
drawing upon and unifying earlier work, re-proved the result 
of Ossona de Mendez~\cite{Osso} on $\alpha$-orientations of a 
graph $G$, mentioning that if $G$ is planar then the lattice of 
$\alpha$-orientations of $G$ 
is isomorphic to the lattice of $\alpha^\perp$-orientations of the dual graph
as defined in my preprint. He then applied his main theorem to Eulerian 
orientations, spanning trees, and Schnyder woods of 3-connected plane graphs.
This line of research was extended by Olivier Bernardi 
and \'Eric Fusy~\cite{BF} in 2011; they showed that the set of Schnyder 
decompositions of a fixed $d$-angulation of girth $d$ has a natural 
distributive lattice structure. 
The 2018 article of Gon\c{c}alves et al.~\cite{GKL} mentioned below 
is also relevant. The most recent and most general presentation 
is the 2025 article of Bernardi, Fusy, and Liang~\cite{BFL}.

In 2005, \'Eric R\'emila presented a fully fleshed-out extension~\cite{Re3}
of the 1992 article of Kenyon and Kenyon~\cite{KK}; working with the Cayley 
graph of a free product of cyclic groups, he gave an algorithm for tiling 
a region using two rectangles of arbitrary size. (A preliminary version of 
this work had appeared in 2002 as a research report coauthored by Bodini 
and R\'emila~\cite{BR}, though only R\'emila appears as author in the 2005 
journal publication.)

In 2007, Knauer~\cite{Knau1} unified the work of Felsner and Propp
and extended the framework to the context of oriented matroids,
though if the oriented matroid is not regular there is not much that 
can be said.  In that same year he published another article~\cite{Knau2}
that generalized the result of Latapy and Magnien~\cite{LaMa}.

In the following year, Knauer teamed up with Felsner to write two papers.
In the first paper~\cite{FK1}, the authors brought height functions
into the broader realm of upper locally distributive lattices,
and more importantly, they unified Knauer's previous work with
the earlier work of Khuller et al.~\cite{KNK} 
through their very general notion of $\Delta$-bonds;
$\Delta$-bonds model planar circular flows, planar $\alpha$-orientations,
$c$-orientations, Schnyder decompositions of a fixed $d$-angulation
of girth $d$, etc. In the second paper~\cite{FK2}, the authors 
further extended their ideas to the context of polytopes.
This work is explained in a more leisurely fashion in~\cite{Knau3}.
$c$-polytopes fall into the family of distributive polytopes
(or $D$-polytopes), where one says that a polytope is distributive
if the coordinatewise max and min of two points in the polytope also 
lie in the polytope, so that the vertex set becomes a distributive lattice.
Furthermore the polytope perspective gives a geometric way to understand
the different finite quotients that can arise from an infinite 
distributive lattice by varying the choice of root vertex or root face;
specifically, it corresponds to intersecting an unbounded polyhedron 
with different subspaces.

In 2008, Olivier Bodini, Thomas Fernique, and \'Eric R\'emila~\cite{BFR}
applied height functions to tilings of infinite regions.

In 2014 Pedro Milet and Nicolau Saldanha~\cite{MS1,MS2,MS3} extended 
earlier work of Saldanha et al.~\cite{STCR} on domino tilings in two dimensions
into three dimensions. Now there are two natural elementary moves 
(flips and trits) that might be used to convert one tiling into another, 
and two topological invariants (flux and twist) that are conserved by moves; 
flux comes from the first homology group and is analogous to the discrete 
divergence mentioned earlier, whereas twist is a discrete analogue of the 
notion of helicity arising in fluid mechanics and topological hydrodynamics.
(See also Milet's thesis~\cite{Mile}.) 

The 2016 article of Pak, Sheffer, and Tassy~\cite{PST}
presented an algorithm for determining
when a simply-connected region can or cannot be tiled by dominoes.
Like Thurston before them,
the authors used height functions in an essential way,
but with clever ideas they were able to improve 
the asymptotic run-time complexity to something essentially optimal;
their procedure (ignoring log factors) runs in time 
proportional to the perimeter of the region rather than the area.
Their procedure works for other sorts of tiling problems as well,
such as lozenge tilings, but they present the method in the case of dominoes
for simplicity of description and ease of analysis. 
The first paragraph of section 5.3 of~\cite{PST} is worthy of note 
for its sketch of the history of height functions;
it discusses such nice articles as Korn and Pak's 2004 
article on T-tetromino tilings~\cite{KP}.

The 2018 article of Daniel Gon\c{c}alves, Kolja Knauer, and 
Benjamin L\'ev\v{e}que~\cite{GKL} extended earlier work on Schnyder woods
and at the same time extended to surfaces of higher genus.
Their homological approach echoed the earlier work of Saldanha 
et al.~\cite{STCR}, although Gon\c{c}alves et al.~were apparently 
unaware of that paper. Their Theorem 4.7 follows from my work, 
but they give a new proof. A more leisurely discussion can be found 
in L\'ev\v{e}que's habilitation thesis~\cite{Leve}.

The 2019 article of Oswin Aichholzer, Jean Cardinal, Tony Huynh, Kolja Knauer,
Torsten M\"utze, Raphael Steiner, and Birgit Vogtenhuber~\cite{ACHK}
study the problem of finding the number of flips required to
convert one $\alpha$-orientation into another. Quite surprisingly,
this problem is shown to be intractable, in the sense that
even determining whether two specified $\alpha$-orientations 
are connected by no more than two flip-moves is NP complete.
See also the followup 2023 article by Cardinal and Steiner~\cite{CS}.
In contrast, the 1995 article of Saldanha et al.~\cite{STCR}
had shown that for domino tilings, the flip distance between
two tilings is essentially the $L^1$-distance between their height-functions.

In 2021, Juliana Freire, Caroline Klivans, Pedro Milet, and 
Nicolau Saldanha~\cite{FKMS}, 
building on earlier work of the third and fourth authors
on three-dimensional domino tilings,
introduced a new move on tilings (refinement), and proved that 
(up to refinement) two tilings are connected by flips and trits 
if and only if they have the same flux, and that (up to refinement) 
two tilings are connected by flips alone if and only if 
they have the same flux and twist. Klivans and Saldanha~\cite{KS}
went on to study analogous notions in higher dimensions.

\subsection*{Acknowledgments}

I thank Noam Elkies, Sergei Fomin, Curtis Greene, David Gupta, Louis Kauffman, 
Greg Kuperberg, William Jockusch, Sampath Kannan, Danny Soroker, Dan Ullman, 
and Lauren Williams for their helpful suggestions.
I am especially grateful to the knowledgeable, industrious, and 
very patient referees who read multiple versions of the article.
The original research was supported by NSA grant MDA904-92-H-3060.
The final write-up process was supported by a Travel Support for 
Mathematicians grant from the Simons Foundation.



\begin{thebibliography}{99}

\bibitem{ACHK}
O.\ Aichholzer, J.\ Cardinal, T.\ Huynh, K.\ Knauer,
T.\ M\"utze, R.\ Steiner, and B.\ Vogtenhuber,
Flip distances between graph orientations,
{\it Algorithmica} {\bf 83} (2021), 116--143.

\bibitem{Andr} 
G.\ Andrews, {\it The Theory of Partitions}, Addison-Wesley, 1976.

\bibitem{ABF}
G.\ Andrews, R.\ J.\ Baxter, and P.\ J.\ Forrester,
Eight-Vertex SOS model and generalized Rogers-Ramanujan-type identities,
{\em J. Stat. Phys.} {\bf 35} (1984), 193--266.

\bibitem{Beij}
H.\ van Beijeren, 
Exactly solvable model for the roughening transition of a crystal surface,
{\it Phys.\ Rev.\ Lett.} {\bf 38} (1977), 993--996.

\bibitem{BF}
O.\ Bernardi and \'E.\ Fusy,
Schnyder decompositions for regular plane graphs and application to drawing,
{\it Algorithmica} {\bf 62} (2012), 1159--1197.

\bibitem{BFL}
O.\ Bernardi, \'E.\ Fusy, and S.\ Liang, Grand-Schnyder Woods,
{\it Ann. Comb.} {\bf 29} (2025), 273--373.

\bibitem{BH}
H.\ W.\ J.\ Bl\"ote and H.\ J.\ Hilhorst,
Roughening transitions and the zero-temperature triangular 
Ising antiferromagnet, 
{\it J.\ Phys.\ A} {\bf 15} (1982), L631--L637.

\bibitem{BFR}
Olivier Bodini, Thomas Fernique, and \'Eric R\'emila,
A characterization of flip-accessibility for rhombus tilings of the whole plane,
{\em Inform. Comput.} {\bf 206} (2008), 1065--1073.

\bibitem{BR}
Olivier Bodini and \'Eric R\'emila,
Tiling a polygon with two kinds of rectangles,
Research Report, LIP RR-2002-46, 
Laboratoire de l’informatique du parall\'elisme. hal-02101939 (2002).

\bibitem{Breh}
E.\ Brehm,
3-orientations and Schnyder 3-tree-decompositions: 
construction and order structure (Diploma thesis), 2000.

\bibitem{BS}
R.\ A.\ Brualdi and M.\ W.\ Schroeder,
Alternating sign matrices and their Bruhat order,
{\it Discrete Math.} {\bf 340} (2017), 1996--2019.

\bibitem{BuPe} 
R.\ Burton and R.\ Pemantle,
Local characteristics, entropy and limit theorems
for uniform spanning trees and domino tilings via transfer impedances,
{\it Ann.\ Probab.} {\bf 21} (1993), 1329--1371.

\bibitem{CS}
J.\ Cardinal and R.\ Steiner,
Inapproximability of shortest paths on perfect matching polytopes,
{\it Math. Program.} Series B (2023).

\bibitem{Ch}
T.\ Chaboud,
Domino tiling in planar graphs with regular and bipartite dual,
{\em Theor. Comput. Sci.} {\bf 159} (1996), 137--142.

\bibitem{CT}
M.\ Cohen and M.\ Teicher,
Kauffman’s clock lattice as a graph of perfect matchings: 
a formula for its height,
{\it Electron. J. Combin.} {\bf 21} (2014), P4.31.

\bibitem{CoLa}
J.\ H.\ Conway and J.\ C.\ Lagarias,
Tilings with polyominoes and combinatorial group theory,
{\it J.\ Combin.\ Theory Ser. A} {\bf 53} (1990), 183--208.

\bibitem{EKLP}
N.\ Elkies, G.\ Kuperberg, M.\ Larsen, and J.\ Propp,
Alternating-sign matrices and domino tilings,
{\it J.\ Algebr.\ Comb.} {\bf 1} (1992), 111--132, 219--234.

\bibitem{Else}
V.\ Elser,
Solution of the dimer problem on a hexagonal lattice with boundary,
{\it J.\ Phys.\ A: Math.\ Gen.} {\bf 17} (1984), 1509--1513.

\bibitem{Felsner}
S.\ Felsner,
Lattice Structures from Planar Graphs,
{\it Electronic J. Combin.} {\bf 11} (2004), \#R15.

\bibitem{FH}
S.\ Felsner and D.\ Heldt,
Mixing Times of Markov Chains on Degree Constrained Orientations 
of Planar Graphs,
{\it Discrete Math. Theor. Comput. Sci.} {\bf 18} (2016), \#20.

\bibitem{FK1}
S.\ Felsner and K.\ B.\ Knauer,
ULD-Lattices and $\Delta$-Bonds,
{\it Combin. Probab. Comput.} {\bf 18} (2009).

\bibitem{FK2}
S.\ Felsner and K.\ B.\ Knauer,
Distributive Lattices, Polyhedra, and Generalized Flow,
{\it European J. Combin.} {\bf 32} (2011), 45--59.

\bibitem{FM}
H.\ de Fraysseix and P.\ Ossona de Mendez,
On topological aspects of orientations,
{\it Discrete Math.} {\bf 229} (2001), 57--72.

\bibitem{FKMS}
J.\ Freire, C.\ Klivans, P.\ H.\ Milet, and N.\ Saldanha,
On the connectivity of spaces of three-dimensional domino tilings,
{\it Trans.\ Amer.\ Math.\ Soc.} {\bf 375} (2022), 1579--1605.

\bibitem{GL}
P.\ Gilmer and R.\ Litherland,
The duality conjecture in formal knot theory,
{\it Osaka J.\ Math.} {\bf 23} (1986), 229--247.

\bibitem{GKL}
D.\ Gon\c{c}alves, K.\ Knauer, and B.\ L\'ev\v{e}que,
On the structure of Schnyder woods on orientable surfaces,
{\it J. Comput. Geom.} {\bf 10} (2019).

\bibitem{Gree}
C.\ Greene, private communication.

\bibitem{GCZ}
D.\ Grensing, I.\ Carlsen, and H.\ Zapp,
Some exact results for the dimer problem 
on plane lattices with non-standard boundaries,
{\it Philos. Mag. A} 
{\bf 41} (1980), 777--781.

\bibitem{HMNS}
M.\ Habiba, R.\ Medina, L.\ Nourinea, and G.\ Steiner,
Efficient algorithms on distributive lattices,
{\em Discrete Appl. Math.} {\bf 110} (2001), 169--187.

\bibitem{HK}
C.\ Hine and T.\ K\'alm\'an,
Clock theorems for triangulated surfaces; preprint (2018), 
available at \arxiv{1808.06091}.

\bibitem{Kauf}
L. Kauffman, {\it Formal Knot Theory},
{\it Mathematical Notes} vol.\ 30, Princeton University Press, 1983.

\bibitem{KK}
C. Kenyon, R. Kenyon, 
Tiling a polygon with rectangles, 
{\em Proc. 33rd Symp. Foundations of Computer Science} (1992), 610--619.

\bibitem{KNK}
S.\ Khuller, J.\ Naor, and P.\ Klein,
The lattice structure of flow in planar graphs,
{\it SIAM J.\ Discrete Math.} {\bf 6} (1993), 477--490.

\bibitem{KS}
C.\ Klivans and N.\ Saldanha,
Domino tilings and flips in dimensions 4 and higher,
{\it Algebr. Comb.} {\bf 5} (2022), 163--185.

\bibitem{Knau1}
K.\ B.\ Knauer,
Partial orders on orientations via cycle flips (master's thesis),
Fachbereich Mathematik der Technischen Universit\"at Berlin, 2007.

\bibitem{Knau2}
K.\ B.\ Knauer,
Distributive lattices on graph orientations, preprint, 2007.
Available at \href
{https://www.ub.edu/comb/koljaknauer/Knauer.pdf}
{https://www.ub.edu/comb/koljaknauer/Knauer.pdf}.

\bibitem{Knau3}
K.\ B.\ Knauer,
Lattices and polyhedra from graphs (doctoral thesis), 2010,
Fakult\"at II – Mathematik und Naturwissenschaften 
der Technischen Universit\"at Berlin.

\bibitem{KP}
M.\ Korn and I.\ Pak,
Tilings of rectangles with T-tetrominoes,
{\em Theor. Comp. Sci.} {\bf 319} (2004), 3--27.

\bibitem{Kupe}
G.\ Kuperberg,
Symmetries of Plane Partitions and the Permanent-Determinant Method,
{\it J.\ Comb. Theory Ser. A} {\bf 68} (1994), 115--151.

\bibitem{LaZh}
P.\ C.\ B.\ Lam and H. Zhang,
A Distributive Lattice on the Set of Perfect Matchings 
of a Plane Bipartite Graph,
{\it Order} {\bf 20} (2003), 13--29.

\bibitem{LaMa}
M.\ Latapy and C.\ Magnien,
Coding Distributive Lattices with Edge Firing Games,
{\it Inform. Process. Lett.} {\bf 83} (2002), 125--128.

\bibitem{LM}
N.\ T.\ T.\ Le and D.\ Mathews,
Generalised Kauffman Clock Theorems, preprint (2025); \arxiv{2510.02911}.

\bibitem{Leve}
B.\ L\'ev\v{e}que,
Generalization of Schnyder woods to orientable surfaces and applications,
habilitation thesis, 2016; \arxiv{1702.07589}.

\bibitem{Levi}
L.\ S.\ Levitov, 
Equivalence of the dimer resonating-valence-bond problem 
to the quantum roughening problem, 
{\it Phys.\ Rev.\ Lett.} {\bf 64} (1990), 92--94.

\bibitem{LMN}
J.\ Linde, C.\ Moore, and M.\ G.\ Nordahl,
An $n$-Dimensional Generalization of the Rhombus Tiling,
{\it Discrete Math. Theor. Comput. Sci. Proc. AA}
(2001), 23--42.

\bibitem{Lova}
L.\ Lov{\'a}sz, {\it Combinatorial Problems and Exercises},
North Holland, 1979.

\bibitem{Mile}
P.\ Milet,
Domino tilings of three-dimensional regions,
Ph.D.\ thesis, 2015; \arxiv{1503.04617}.

\bibitem{MS1}
P.\ Milet and N.\ Saldanha,
Domino tilings of three-dimensional regions: flips, trits and twists,
\arxiv{1410.7693}.

\bibitem{MS2}
P.\ Milet and N.\ Saldanha,
Twists for duplex regions, \arxiv{1411.1793}.

\bibitem{MS3}
P.\ Milet and N.\ Saldanha,
Flip invariance for domino tilings of three-dimensional regions with two floors,
{\it Discrete Comput.\ Geom.} {\bf 53} (2015), 914--940.

\bibitem{MRR1}
W.H.\ Mills, D.P.\ Robbins, and H.\ Rumsey, Jr.,
Alternating sign matrices and descending plane partitions,
{\it J.\ Comb.\ Theory Ser. A} {\bf 34} (1983), 340--359.

\bibitem{MN}
C.\ Moore and M.\ E.\ J.\ Newman,
Height Representation, Critical Exponents, and Ergodicity in the Four-State
Triangular Potts Antiferromagnet, 
{\it J.\ Stat.\ Phys.} {\bf 99} (2000), 629--660.

\bibitem{MNMS}
C.\ Moore, M.\ G.\ Nordahl, N.\ Minar, and C.\ R.\ Shalizi,
Vortex dynamics and entropic forces in antiferromagnets 
and antiferromagnetic Potts models,
{\it Phys.\ Rev.\ {\bf E}} {\bf 60} (1999), 5344.

\bibitem{MRR2}
C.\ Moore, I.\ Rapaport and \'E. R\'emila,
Tiling groups for Wang tiles,
{\it Proc. ACM-SIAM Symp. Discrete Algorithms} {\bf 13}
(SODA 2002), 402--411.

\bibitem{Osso}
P.\ Ossona de Mendez,
Orientations bipolaires (Ph.D.\ thesis),
\'Ecole des Hautes \'Etudes en Sciences Sociales, Paris, 1994.

\bibitem{Pak1}
I.\ Pak,
Ribbon Tile Invariants,
{\it Trans. Amer. Math. Soc.} {\bf 352} (2000), 5525--5561.

\bibitem{Pak2}
I.\ Pak,
Tile invariants: new horizons,
{\it Theor. Comput. Sci.} {\bf 303} (2003), 303--331.

\bibitem{PST}
I.\ Pak, A.\ Sheffer, and M.\ Tassy, Fast domino tileability,
{\it Discrete Comput.\ Geom.} {\bf 56} (2016), 377-394.

\bibitem{Pr0}
O.\ Pretzel,
On graphs that can be oriented as diagrams of ordered sets,
{\it Order} {\bf 2} (1985), 25--40.

\bibitem{Pr1}
O.\ Pretzel,
On reorienting graphs by pushing down maximal vertices,
{\it Order} {\bf 3} (1986), 135--153.

\bibitem{Pr2}
O.\ Pretzel,
Orientations and reorientations of graphs,
in ``Combinatorics and Ordered Sets'',
{\it Contemp.\ Math.} {\bf 57} (1986), 103--125.

\bibitem{Pr3}
O.\ Pretzel,
On reorienting graphs by pushing down maximal vertices--II,
{\it Discrete Math.} {\bf 270} (2003), 227--240.

\bibitem{Prop}
J.\ Propp,
The many faces of alternating-sign matrices, 
{\it Discrete Math. Theor. Comput. Sci. Proc. AA}
(2001), 43--58.

\bibitem{PW}
J.\ G.\ Propp and D.\ B.\ Wilson, 
Exact sampling with coupled Markov chains and applications 
to statistical mechanics, 
{\it Random Struct. Algorithms} {\bf 9} (1996), 223--252.

\bibitem{Re1}
\'E.\ R\'emila,
Tiling Groups: New Applications in the Triangular Lattice,
{\it Discrete Comput. Geom.} {\bf 20} (1998), 189--204.

\bibitem{Re2}
\'E.\ R\'emila,
The lattice structure of the set of domino tilings of a polygon, 
{\it Theor. Comp. Sci.} {\bf 322} (2004), 409--422.

\bibitem{Re3}
\'Eric R\'emila,
Tiling a polygon with two kinds of rectangles,
{\it Discrete Comput. Geom.} {\bf 34} (2005), 313-330.

\bibitem{Robb}
D.P.\ Robbins,
The story of 1, 2, 7, 42, 429, 7436, ...,
{\it Math.\ Intelligencer} {\bf 13} (1991), 12--19.

\bibitem{SZ}
H.\ Sachs and H.\ Zernitz,
Remark on the dimer problem,
{\it Discrete Appl. Math.} {\bf 51} (1994), 171--179.

\bibitem{STCR}
N.C.\ Saldanha, C.\ Tomei, M.A.\ Casarin, Jr., and D.\ Romualdo,
Spaces of Domino Tilings,
{\it Discrete Comput.\ Geom.} {\bf 14} (1995), 207--233.

\bibitem{Shef}
S.\ Sheffield,
Ribbon tilings and multidimensional height functions,
{\it Trans. Amer. Math. Soc.} {\bf 354} (2002), 
4789--4813.

\bibitem{Temp}
H.\ N.\ V.\ Temperley, in:
{\it Combinatorics},
London Math.\ Soc.\ Lecture Notes Series {\bf 13} (1974), 202--204.

\bibitem{Thom}
M.\ Thompson, private communication.

\bibitem{Thur}
W.\ Thurston,
Conway's tiling groups,
{\it Amer. Math. Monthly} {\bf 97} (1990), 757--773.

\end{thebibliography}
\end{document}